\documentclass{aupl}

\usepackage{
amsfonts,
latexsym,
amssymb,
amsmath,
amsthm,
enumerate,
verbatim,
mathrsfs,
}

\usepackage{url}

\newcommand{\labbel}[1]{\label{#1} [[{\bf #1}]]}  %
\renewcommand{\labbel}{\label}

\newtheorem{theorem}{Theorem}[section]
\newtheorem{lemma}[theorem]{Lemma}

\newtheorem{proposition}[theorem]{Proposition} 
 
\newtheorem{corollary}[theorem]{Corollary} 
 
\newtheorem*{claim}{Claim} 
\newtheorem*{subclaim}{Subclaim}

\newtheorem*{theorem*}{Theorem}
\newtheorem*{corollary*}{Corollary}

\theoremstyle{definition}
\newtheorem{definition}[theorem]{Definition}

\newtheorem{problem}[theorem]{Problem}

\theoremstyle{remark}
\newtheorem{remark}[theorem]{Remark}
 
\newtheorem{example}[theorem]{Example}

\newtheorem{observation}[theorem]{Observation} 
\newtheorem{notation}[theorem]{Notation}

\numberwithin{equation}{section}

\newcommand{\brfrt}{\hspace{0 pt}}

\newcommand{\hyphe}{\frac{\ }{\ }}

\newcommand{\timee}{{\hspace {0.4pt} \times \hspace {0.4pt}}}

 \allowdisplaybreaks[0]
\begin{document}

\title{Mitschke's Theorem is sharp}

\author{Paolo Lipparini} 
\address{Dipartimento 
 di Matematica\\Viale della  Ricerca  Quasi Unanime
\\Universit\`a di Roma ``Tor Vergata'' 
\\I-00133 ROME ITALY}
\urladdr{http://www.mat.uniroma2.it/\textasciitilde lipparin}


\thanks{Work performed under the auspices of G.N.S.A.G.A. Work 
partially supported by PRIN 2012 ``Logica, Modelli e Insiemi''.
The author acknowledges the MIUR Department Project awarded to the
Department of Mathematics, University of Rome Tor Vergata, CUP
E83C18000100006.
\\
\indent Date: \today}

\begin{abstract}
A.\ Mitschke
showed that a variety with an 
$m$-ary near-unanimity 
term has J{\'o}nsson terms $t_0, \dots,  t _{2m-4} $
witnessing congruence distributivity. 
We show that Mitschke's result is sharp.
 We also evaluate the best possible number of Day terms 
witnessing congruence modularity.
More generally, we characterize exactly
the best bounds for  many
congruence identities satisfied by varieties 
with an $m$-ary near-unanimity 
term.
\end{abstract}

\maketitle

\section{Introduction} \labbel{intro}

Recall that a term $u$ is a \emph{near-unanimity term}  
(in some algebra or in some variety) if  all the equations of the form
\begin{equation*} 
u(x,x, \dots, x,y,x, \dots, x,  x) =x
\end{equation*}    
are satisfied, with just one occurrence of $y$
in any possible position. 

Near-unanimity terms have been originally studied in connection
with generalizations of the Chinese remainder theorem \cite{BP}.
More recent research has shown connections with computational
complexity, e.~g. \cite{Ba,BK,BIM}.
Joins of varieties with near-unanimity terms
have been studied in 
\cite{CCV}. 

 A.\ Mitschke \cite{Mi}
proved that every variety $\mathcal {V}$ with a near-unanimity term 
is congruence distributive.
A more direct proof, credited to E. Fried,
can be found in  Kaarli
and Pixley \cite[Lemma 1.2.12]{KP}. 
Compare also Barto and Kozik \cite[Section 5.3.1]{BK} 
and \cite[Section 5]{B}. 

In particular, any  variety with a near-unanimity term  
is congruence modular.
The distributivity \cite{JD}  and 
modularity \cite{D} levels of varieties with a near-unanimity term 
have been evaluated.

\begin{theorem} \labbel{mnun}
Let $m \geq 3$.

  \begin{enumerate}    
\item    
(Mitschke \cite{Mi})
 A variety with an
$m$-ary near-unanimity term is $2m{-}4$-distributive.
\item
(Sequeira \cite[Theorem 3.19]{S})
A variety with an
$m$-ary near-unanimity term is $2m{-}3$-modular.
 \end{enumerate}
 \end{theorem}

In this note we   show that Theorem \ref{mnun}
gives the best possible evaluations.
Section \ref{mainc} presents  our main construction, where
we build  appropriate subalgebras of certain products.
 The construction is then iterated in 
Section \ref{Mitschkesharp} in order to get 
counterexamples showing that 
Theorem \ref{mnun} cannot  be improved.
In Section \ref{expl} 
we exemplify our methods by
presenting some more concrete examples. 
Section \ref{expl} is largely self-contained.
Further remarks are contained in Section \ref{fur}.

We shall assume familiarity with the basic notions
of universal algebras, as presented, e.~g., in
\cite{MMT}. 
The notions we shall use
admit equivalent reformulations in terms of congruence identities,
as given by the following table.
\begin{equation}\labbel{cacu}     
\begin{aligned}
& \text{ $n$-distributive  \quad\quad }   
&&
 \alpha ( \beta \circ \gamma )
 \subseteq  
\alpha \beta \circ \alpha \gamma \circ {\stackrel{n}{\dots}} 
\\ 
& \text{ $n$-alvin  }   
&& 
\alpha ( \beta \circ \gamma )
 \subseteq 
 \alpha \gamma  \circ  \alpha \beta  \circ {\stackrel{n}{\dots}} 
\\
& \text{ $n$-modular  }   
&&
 \alpha ( \beta \circ \alpha  \gamma \circ \beta )
 \subseteq  
\alpha \beta \circ \alpha \gamma \circ {\stackrel{n}{\dots}} 
\\
& \text{ $n$-reversed-modular \quad \quad }   
&&
 \alpha ( \beta \circ \alpha  \gamma \circ \beta )
 \subseteq  
\alpha \gamma  \circ \alpha \beta  \circ {\stackrel{n}{\dots}} 
\end{aligned}
  \end{equation}    
A notion mentioned on the left holds 
in some  variety $\mathcal {V}$ 
if and only if $\mathcal {V}$ satisfies the corresponding 
congruence identity on the right, that is, 
the identity  holds for every algebra $\mathbf A$ in $\mathcal {V}$ 
and all congruences in $\mathbf A$.
In the above formulae juxtaposition
denotes intersection.
For $\varepsilon$ and $\delta$ binary relations,  
$\varepsilon \circ \delta \circ {\stackrel{k}{\dots}}  $
denotes the relation $ \varepsilon \circ \delta \circ \varepsilon \circ \delta \circ \dots$
with $k$ factors, that is, $k-1$ occurrences of  $ \circ $.
If, say, $k$ is even, then we  write
$\varepsilon \circ \delta \circ {\stackrel{k}{\dots}} \circ \delta  $
when we want to make clear that $\delta$ is the last factor.
Sometimes, for readability or convenience, 
 we might add further factors in the above expressions,
as in  $\varepsilon \circ \delta \circ \varepsilon \circ {\stackrel{k}{\dots}}
\circ \varepsilon \circ \delta   $. In any case, the number above the dots
represents the number of 
occurrences of  $ \circ $ minus one.

Usually, the notions introduced in
\eqref{cacu} are defined in an equivalent way (throughout a variety)
by means of the existence of a certain number of terms,
called after B. J{\'o}nsson, for $n$-distributive, and
after A. Day, for $m$-modular. The terms for alvin and reversed modularity
are obtained from J{\'o}nsson and Day terms, respectively,
by exchanging the conditions 
for even and odd indices. As in \cite{day}, such 
``reversed'' conditions will be the key to the exact
evaluation of the appropriate parameters: 
 dealing
only with the more standard notions
we could not succeed in proving the exact results.  
In any case, here we shall not need terms for distributivity and modularity;
dealing with congruence identities will prove much simpler. 
The reasons why in certain cases congruence identities or even relation
identities are  more convenient than terms are explained
in Tschantz \cite{Ts}.
 
It is immediate from \eqref{cacu}
that every $n$-alvin variety is $n{+}1$-distributive, that
every $n$-distributive variety is $n{+}1$-alvin,
and corresponding results hold for modularity and
reversed modularity.

See \cite{day}, in particular, Section 2 therein,
 for a full discussion of the equivalences
presented in table \eqref{cacu}. 
See   \cite{B,day} for more results related or similar to Theorem \ref{mnun}
and for further comments.
In particular in \cite{day}, among many other examples,
 we  presented, for every even $n \geq 2$, the example
of a locally finite $n$-distributive variety which is neither  
$2n{-}2$-modular, nor  $2n{-}1$-reversed-modular.

\section{The main construction} \labbel{mainc}

\begin{definition} \labbel{kmaj}
Let $m \geq 3$ and $1\leq k \leq m $.

If some algebra $\mathbf A$ has a special element
$0$,  we say that $0$ is a \emph{$k$-absorbing
element for a  term $u$} if  $u(a_1, \dots, a_m) = 0$,
whenever $0$  occurs at least $k$-times
in the arguments of $u$, more formally, whenever
$|\{ i \mid a_i=0 \}| \geq k$.

A term $u$ is 
a \emph{$k$-majority term} in some algebra $\mathbf A$ 
(in some variety $\mathcal {V}$)
if every element of $\mathbf A$ 
(of every algebra in $\mathcal {V}$)
is $k$-absorbing
for $u$.   
In other words, a $k$-majority term
is supposed to satisfy the equation
 $u(x_1, \dots, x_m) = x$,
whenever the variable $x$   occurs at least $k$-times
in the arguments of $u$.
Clearly, $k> \frac{m}{2} $,
unless we are in a trivial variety. 

An $m$-ary term $u$ 
is \emph{idempotent} if it is 
an $m$-majority term, namely, if
the equation $u(x, x, \dots, x) = x$ is satisfied.
 
An  $m$-ary term is a \emph{near-unanimity term} 
if it is  
an $m-1$-majority term. 
 
An $m$-ary term $u$ is 
\emph{symmetrical} in some algebra $\mathbf A$ 
(in some variety $\mathcal {V}$)
if all the equations $u(x_1, \dots, x_m)=
 u(x_{ \tau  (1)}, \dots, x_ { \tau  (m)})$
hold  in  $\mathbf A$ 
 (in $\mathcal {V}$),
for all permutations $ \tau $ of
$\{ 1, \dots, m\}$. 
\end{definition}    

In principle, when $k< m-1$, the notion of a    
 $k$-majority $m$-ary term has little interest,
since it implies the existence of a near-unanimity term 
of arity $<m$.   However, we shall merge
different varieties with a $k$-majority term, for distinct values of $k$,
in such a way that the resulting variety  $\mathcal {V}$ 
has an  $m{-}1$-majority term 
(namely, a near-unanimity term)
and provides all the desired counterexamples.

The next construction and, more generally, all the arguments in the present note
share many aspects in common with the constructions we have  performed
in \cite{day}.
However, an important difference should be mentioned. 
In the  constructions in \cite{day}, at each inductive step,
we have taken the product of some formerly constructed algebra $\mathbf A_4$ 
with three further algebras. One of these additional algebras,
the algebra $\mathbf A_3$ in \cite{day}, 
is a term-reduct of the two-elements lattice $\mathbf C_2$.
Thus in \cite{day} at each induction step a reduct of
$\mathbf C_2$ is added as a new factor.
In the present situation, instead, it is necessary to fix 
 the $\mathbf C_2$-reduct
once and for all at the beginning, hence here the induction steps
start with  a subalgebra of $\mathbf A_3 \times\mathbf A_4$.
Let us also mention that, for convenience,
 here we shall shift the third and fourth indices,
in comparison with \cite{day}. In particular, the reduct of  
$\mathbf C_2$ here will appear at the fourth place.

We shall frequently consider special elements
$0_z\in \mathbf A_z$, for $z=1,2,4$.
When no confusion is possible, 
we shall omit the subscripts.
The types 
introduced in the next lemma 
have been used also  in many constructions
from \cite{day}. 
Since, as we mentioned,
we are shifting the last two coordinates,
the correspondence with \cite{day}
is exact only modulo a permutation 
of the coordinates. This is the reason why
the types here are denoted by, say, II$^ \sigma $,
rather than II.

Throughout the present note 
objects like $\mathbf A_1 \times \mathbf A_2 \times 
\mathbf A_3 \times \mathbf A_4$
and, say, $\mathbf A_1 \times \mathbf A_2 \times 
(\mathbf A_3 \times \mathbf A_4)$
shall be always identified, namely, we  consider them 
modulo isomorphism through the natural correspondence.

\begin{lemma} \labbel{lemnu}
Suppose that $\mathbf A_1$, $\mathbf A_2$, $\mathbf A_3$
and $\mathbf A_4$ are algebras with exactly
one $m$-ary operation $u$. Suppose that $ 3 \leq m $, 
$1 \leq h \leq k $ and
$h+k \leq m$. Suppose further that
$0_z\in \mathbf A_z$, for $z=1,2,4$  and
  \begin{enumerate}   
 \item  
$0_z $ is $h$-absorbing for $u$ in $  \mathbf A_z$, for $z=1,2$, 
 \item
$u$ is a $k$-majority term in $\mathbf A_3$, and
\item  
$0_4 $ is $2$-absorbing for $u$ in $  \mathbf A_4$. 
 \end{enumerate} 

Suppose that $a, d \in  A_3$,
$\mathbf F$ is a subalgebra of  $\mathbf A_3 \times \mathbf A_4$ 
and let $B=B(a,d)$ be the subset of $  A_1 \times A_2 \times F $
consisting of the elements which have one of the following types
\begin{equation*}
\begin{gathered}
\text{Type I$^ \sigma$} \\
(\hyphe, 0, a, \hyphe)
\end{gathered} 
\qquad\qquad
\begin{gathered}
\text{Type II$^ \sigma$} \\
(0, 0, \hyphe, \hyphe),
  \end{gathered} 
\qquad\qquad
\begin{gathered}
\text{Type III$^ \sigma$} \\
(0, \hyphe, d, \hyphe)
  \end{gathered} 
\qquad\qquad
\begin{gathered}
\text{Type IV$^ \sigma$} \\
(\hyphe, \hyphe,  \hyphe, 0),
  \end{gathered} 
\end{equation*}    
where dashed places
can be filled with arbitrary elements from the 
corresponding algebras,
and under the provision that each $4$-uple 
actually belongs to $  A_1 \times A_2 \times F $,
namely, that the couple consisting of the last two coordinates
belongs to  $F$.
Recall that  
 we are  omitting the subscripts relative to the  $0$'s.

 Then $B=B(a,d)$ is the base set for a subalgebra
$\mathbf B = \mathbf B(a,d)$  of 
$\mathbf A_1 \times \mathbf A_2 \times \mathbf F$,
hence also a subalgebra  of $\mathbf A_1 \times \mathbf A_2 \times 
\mathbf A_3 \times \mathbf A_4$.
 \end{lemma}

 \begin{proof} 
First notice that $B$ is nonempty,
since there exists at least an element of type II$^ \sigma$. 
Suppose that $b_1, \dots, b_m \in B$.
We have to show that
 $b=u(b_1, \dots, b_m) \in B$.
Since 
$\mathbf A_1 $ and $  \mathbf A_2$ are algebras and
$\mathbf F$ is a subalgebra of  $\mathbf A_3 \times \mathbf A_4$,
if each $b_i$ belongs to  
$A_1 \times  A_2 \times F$,
then $b \in A_1 \times  A_2 \times F$.
Hence it remains to show that $b$
has one of types I$^ \sigma$ - IV$^ \sigma$. 

If at least two $b_i$'s have type IV$^ \sigma$, then
$b$ has type IV$^ \sigma$, by (3), 
hence we can suppose that at most one   
$b_i$ has type IV$^ \sigma$.

If at least $h$-many $b_i$'s have type I$^ \sigma$ or II$^ \sigma$
(hence have $0$ in the second position) 
and at least $h$-many $b_i$'s have type III$^ \sigma$ or II$^ \sigma$
(hence have $0$ in the first position),
then
$b$ has type II$^ \sigma$, by (1), and we are done in this case, as well.

Otherwise, 
there are, say,
at most $h{-}1$-many $b_i$'s having type
either  I$^ \sigma$ or II$^ \sigma$.
Since we have assumed that at most 
one $b_i$ has type IV$^ \sigma$,
then there are at least 
$m-(h-1)-1 = m-h \geq k$  many
$b_i$'s of type III$^ \sigma$.
Then
$b$ has type III$^ \sigma$, by 
(1), $h \leq k$ and (2).
Symmetrically,
if there are 
at most $h-1$-many $b_i$'s of type III$^ \sigma$ or II$^ \sigma$,
then there are at least 
$k$-many
$b_i$'s of type I$^ \sigma$,
thus
$b$ has type I$^ \sigma$.
\end{proof}

\section{Mitschke's Theorem is sharp} \labbel{Mitschkesharp} 

For fixed $m \geq 3$, we now consider lattice terms of the form
\begin{equation*}     
u _{2,m} = \prod _{i < j<m}  (x_i + x_j) , \qquad
u _{3,m} = \prod _{ i<j<k <m} (x_i + x_j + x_k) 
 \end{equation*}
and so on.
We shall combine various lattice reducts
defined using the above terms
in order to obtain our counterexamples.
 As shown by the following remarks, it is 
not enough to consider just one of the above term-reducts.

\begin{example} \labbel{ex}
For  $m \geq 3$,  consider the  term-reduct $\mathcal V_m^d$ of the 
variety of distributive lattices with the only 
$m$-ary operation corresponding to the term 
$u _{2,m}$ defined above.
Trivially the operation in $\mathcal V_m$ is 
an $m$-ary near-unanimity term, since 
  $u _{2,m}$ is an $m$-ary near-unanimity term
in lattices.  Moreover, $\mathcal V_m^d$ is locally finite, 
being a term-reduct of a locally finite variety.

It is easy to see that if $m\geq 4$,
then $\mathcal V_m^d$ has not an $m{-}1$-ary
near-unanimity term, thus, in general,
the existence of an $m$-ary
near-unanimity term does not imply an $m{-}1$-ary
near-unanimity term, even for locally finite varieties.
See \cite[Lemma 3.4]{S} for a slightly more involved example
(not locally finite).

To check that
$\mathcal V_m^d$ has not an $m{-}1$-ary
near-unanimity term, let $\mathbf A $
be the $u _{2,m}$ term-reduct of  
$  \mathbf C_2 \times \mathbf C_2 \times \dots \times \mathbf C_2$
with $m-1$ factors, 
where  $\mathbf C_2$ is the two-elements lattice
with universe $\{ 0,1\}$. 
Then $B=A \setminus (1,1, \dots, 1)$
is the universe for a subalgebra of $\mathbf A$,
since, for any $m$-uple  of elements of $B$,
we have at least two elements with a $0$  
at the same component, hence we still get $0$
at that component when applying $u _{2,m}$. 
On the other hand, $\mathcal V_m^d$  has not
an $m{-}1$-ary near-unanimity term, since,
were $v$  such a term, then in $\mathbf A$ 
\begin{align*}
&v((0,1,1, \dots, 1), (1,0,1, \dots, 1), \dots, (1,1,1, \dots, 0))=
\\
&(v(0,1,1, \dots, 1), v(1,0,1, \dots, 1), \dots, v(1,1,1, \dots, 0))=
(1,1,1,\dots,1),
  \end{align*}    
contradicting the above-proved fact that $\mathbf  B$ is a 
subalgebra of $\mathbf A$, since 
$(1,1,1, \allowbreak \dots,1) \notin B$.

If $\mathcal V_m$ is the corresponding term-reduct of the variety 
of all lattices, then $\mathcal V_m$, too,
has   an $m$-ary near-unanimity term 
but not an $m{-}1$-ary near-unanimity term.
In this case,  $\mathcal V_m$ is not locally finite. 
\end{example}

\begin{remark} \labbel{rmb}   
In general, we cannot use the above example
in order to show that Theorem \ref{mnun} 
is the best possible result.
Indeed, Baker \cite{Bk} showed that any congruence
distributive term reduct of lattices is $4$-distributive.
While $\mathcal V_4$ is actually 
not $3$-distributive and not $4$-modular \cite{B}, thus 
$\mathcal V_4$ indeed shows that 
Theorem \ref{mnun} gives the best possible result for $m=4$,
Baker's Theorem prevents $\mathcal V_m$
to be a suitable counterexample for larger $m$.

Henceforth  a more involved approach is necessary.  

Notice that  
Baker's result can be generalized to the effect that 
any congruence
distributive term reduct of Boolean algebras is $4$-distributive.
This statement is immediate 
from  \cite[Theorem 6.4(3)]{B}.
Hence, for our present purposes, considering Boolean algebras
in place of lattices provides no special advantage.
In this respect, however, see Proposition \ref{MM}.  
\end{remark}

\begin{definition} \labbel{latterms}
Suppose that  $m \geq 3$ and $1 \leq j \leq m$.

Let $u _{j,m}$
be  the following 
$m$-ary 
lattice term  
\begin{equation}\labbel{lte}
u _{j,m}(x_1, \dots, x_m)= \prod _{|J|=j} \sum _{i \in J} x_i 
   \end{equation}    
where $J$ varies on  subsets of 
$\{ 1, \dots, m \}$.

Of course, 
strictly speaking,  $u _{j,m}$, as a term, is uniquely defined only modulo some
fixed arrangement of summands and factors.
 However, we shall be only interested on $u _{j,m}$
as an operation, hence the actual syntactical definition of the term
$u _{j,m}$ shall not be relevant in what follows. 
\end{definition}
   
\begin{observation} \labbel{obs} 
Notice that, in every lattice with minimum,
the minimum  $0$ is $j$-absorbing for $u _{j,m}$.
Moreover, in every lattice, 
$u _{j,m}$ is a $p$-majority term
for $p= \max \{j,  m{-}j{+}1\}$.
In particular, if
$j \leq \frac{m+1}{2} $,
then  $u _{j,m}$ is an $ m{-}j{+}1$-majority term.
 \end{observation}   

\begin{definition} \labbel{vars}
Suppose that  $m \geq 3$ and $2 \leq j < m$.  
If $\mathbf L$ is a lattice, let 
$\mathbf L ^{{\rm nu}, j,m} $ be the term-reduct of 
  $\mathbf L$ with $u _{j,m}$ as the only operation
(henceforth always named as $u$).
Let $\mathbf N ^{j,m} = \mathbf C_2 ^{{\rm nu},j,m} $,
where $\mathbf C_2$ is the two-elements lattice
with base set $\{ 0,1\}$. 

Let 
$\ell = \frac{m+1}{2} $
if $m$ is odd, and
  $\ell = \frac{m}{2} $
if $m$ is even. Let
$\mathcal {N}_m$ be the variety 
generated by the algebras
\begin{equation*}
 \mathbf N ^{2, m}, \quad 
\mathbf N ^{3, m}, \quad 
 \dots, \quad \mathbf N ^{ \ell, m}.
\end{equation*}

The definition is well-posed
since the second superscript  
determines the type of the algebra, 
in the present case, the arity of the only operation.
 \end{definition}   

Conventionally, we let 
$ \varepsilon  \circ   \delta  \circ {\stackrel{1}{\dots}}
= \varepsilon   $ and
$  \varepsilon  \circ  \delta  \circ {\stackrel{0}{\dots}}
= 0$, where $0$ is the minimal congruence in the algebra under consideration.
If $R$ is a binary relation,
$R^k$ denotes  $ R \circ   R  \circ {\stackrel{k}{\dots}} \circ R$.

\begin{theorem} \labbel{mmm}
Let $m \geq 3$.
The variety $\mathcal {N}_m$ 
is locally finite and has
an $m$-ary symmetrical near-unanimity term. Moreover 
  \begin{enumerate}    
\item    
$\mathcal {N}_m$ is not $2m{-}4$-alvin,
in particular, not $2m{-}5$-distributive.
\item
More generally, the following congruence identity
fails in  $\mathcal {N}_m$
\begin{equation}\labbel{cid}
\alpha (\beta \circ \gamma ) 
\subseteq
(\alpha ( \gamma \circ \beta )) ^{m-2} .  
   \end{equation}    
\item
$\mathcal {N}_m$ is not $2m{-}3$-reversed-modular,
in particular, not $2m{-}4$-modular.
\item
Still more generally, the following congruence identity
fails in  $\mathcal {N}_m$
\begin{equation}\labbel{cidd}
\alpha (\beta \circ ( \alpha \gamma \circ \alpha \beta 
\circ {\stackrel{q-2}{\dots}} \circ \alpha \beta ^\bullet   )
\circ   \gamma ^ \bullet) 
\subseteq
(\alpha ( \gamma \circ \beta \circ {\stackrel{q}{\dots}}
 \circ  \beta ^ \bullet  )) ^{m-2},  
   \end{equation}    
for every $q \geq 2$,
 where 
$\beta^ \bullet = \beta $,
$ \gamma ^ \bullet = \gamma  $
if $q$ is even and 
$\beta^ \bullet = \gamma $,
$ \gamma ^ \bullet = \beta   $
if $q$ is odd.
 \end{enumerate}
 \end{theorem}   

\begin{proof}
The variety $\mathcal {N}_m$ has
an $m$-ary near-unanimity term, actually, 
an $m$-ary near-unanimity operation,
since 
in each algebra 
$ \mathbf N ^{j, m} $,
for $2 \leq j \leq \ell$,  
the only operation is a near-unanimity operation. 
Indeed, by Observation \ref{obs},
the operation of  $ \mathbf N ^{j, m} $ is a
$p$-majority term
for $p= \max \{j,  m{-}j{+}1\}$.
Since  $2 \leq j \leq \ell \leq \frac{m+1}{2} $
and $m \geq 3$, 
we have $p \leq m-1$, for every $j$ in the interval under 
consideration.
Now notice that if $p\leq p'$,
then a  $p$-majority term
is a $p'$-majority term.
Hence  in each algebra 
$ \mathbf N ^{j, m} $
the operation is an $m{-}1$-majority term,
that is, a near-unanimity term.
The operation is symmetrical, since it is symmetrical on
each generating algebra.

Since the variety of distributive lattices is locally finite, 
each algebra $ \mathbf N ^{j, m}$ 
generates a locally finite variety, hence 
$\mathcal {N}_m$ is locally finite, being the join
of a finite number of locally finite varieties.

We now show that (1) - (3)
all follow from (4).
Of course, the reader interested only in (1) - (3)
might work out the details of the following arguments
in the corresponding simplified setting. Cf. also Section \ref{expl} below.

To show that (2) follows from (4) observe that   
 \eqref{cid} is the special case 
$q=2$ of \eqref{cidd}.
Moreover, (1) is immediate from (2), since
$ \alpha \beta \circ \alpha \gamma \subseteq \alpha (\beta \circ \gamma )   $.
Recall the  conditions given by \eqref{cacu}. 
To show that 
(4) implies (3), we first establish 
a condition of independent interest and which, for $q $ odd,
is equivalent to 
\eqref{cidd} in every algebra.  

\begin{lemma} \labbel{cilem}
If $m,q \geq 3$ and $q$ is odd
then identity \eqref{cidd} is equivalent to
\begin{equation}\labbel{ciodd}
\alpha (\beta \circ ( \alpha \gamma \circ \alpha \beta 
\circ {\stackrel{q-2}{\dots}} \circ \alpha \gamma    )
\circ   \beta ) 
\subseteq
\alpha \gamma  \circ
\big(\alpha (  \beta \circ \alpha \gamma \circ  
\beta \circ {\stackrel{q-2}{\dots}} \circ 
 \alpha \gamma \circ  \beta ) \circ \alpha \gamma \big) ^{m-2}  
   \end{equation}    
in every algebra. 
 \end{lemma}

Indeed, 
by taking $\alpha \gamma $ in place of $\gamma$ in 
\eqref{cidd} and since $q$ is odd, we get
\begin{align*} 
&\alpha (\beta \circ ( \alpha \gamma \circ \alpha \beta 
\circ {\stackrel{q-2}{\dots}} \circ \alpha \gamma    )
\circ   \beta ) 
\subseteq
\big(\alpha ( \alpha \gamma \circ \beta \circ \alpha \gamma 
 \circ {\stackrel{q}{\dots}}
 \circ \beta \circ   \alpha \gamma   )\big) ^{m-2}=
\\
&\alpha \gamma \circ \alpha (\beta \circ \alpha \gamma 
 \circ {\stackrel{q-2}{\dots}} \circ \beta) \circ \alpha \gamma \circ 
\alpha \gamma \circ \alpha (\beta \circ \alpha \gamma 
 \circ {\stackrel{q-2}{\dots}} \circ \beta) \circ \alpha \gamma  \dots
\\
&\dots \alpha \gamma \circ \alpha (\beta \circ \alpha \gamma 
 \circ {\stackrel{q-2}{\dots}} \circ \beta) \circ \alpha \gamma
\circ \alpha \gamma 
\circ \alpha (\beta \circ \alpha \gamma 
 \circ {\stackrel{q-2}{\dots}} \circ \beta) \circ \alpha \gamma =
\\
&\alpha \gamma  \circ
\big(\alpha (\beta \circ \alpha \gamma 
 \circ {\stackrel{q-2}{\dots}} \circ \beta)
 \circ \alpha \gamma \big) ^{m-2}  
\end{align*}   
since    $\alpha ( \alpha \gamma \circ \beta \circ \alpha \gamma 
 \circ {\stackrel{q}{\dots}}
 \circ \beta \circ   \alpha \gamma   ) = \alpha \gamma \circ \alpha (\beta \circ \alpha \gamma 
 \circ {\stackrel{q-2}{\dots}} \circ \beta) \circ \alpha \gamma$
and $ \alpha \gamma \circ \alpha \gamma = \alpha \gamma $,
both   $\alpha$ and $\alpha \gamma $ being equivalence relations. 
Hence \eqref{cidd} implies \eqref{ciodd}. 

On the other hand, for all congruences $\alpha$, $\beta$ and $\gamma$,
\begin{equation*}     
\alpha \gamma  \circ
\alpha (  \beta \circ \alpha \gamma \circ  
\beta \circ {\stackrel{q-2}{\dots}} \circ
 \alpha \gamma \circ  \beta ) \circ \alpha \gamma 
\subseteq 
\alpha (\gamma  
\circ  \beta \circ \gamma \circ  
\beta \circ {\stackrel{q}{\dots}} \circ
 \gamma \circ  \beta  \circ  \gamma ),
 \end{equation*} 
thus
\eqref{ciodd} implies \eqref{cidd}, for $q$ odd.    
We have proved Lemma \ref{cilem}.

If $q=3$, then $q-2=1$,
hence in this case \eqref{ciodd}
becomes exactly the condition \eqref{cacu}  for 
$2m{-}3$-reversed modularity,
thus, by Lemma \ref{cilem},
the special case $q=3$ of clause
(4) in Theorem \ref{mmm} implies clause (3).

Since we have showed that clause (4)
implies all the other clauses in \ref{mmm}, 
it remains to prove (4).
The proof shall involve
further definitions, notation, claims
and shall go through a finite induction
divided in three steps,  using 
Lemma \ref{lemnu}. 

Fix $m \geq 3$, $q \geq 2$ 
and let $\ell$ be as in Definition \ref{vars}. 
For every $j$ with  $2 \leq j \leq \ell$,
let $\mathcal {N}_m^j$ be the variety 
generated by the algebras
\begin{equation*}
 \mathbf N ^{j, m}, \quad 
\mathbf N ^{j+1, m}, \quad 
 \dots, \quad \mathbf N ^{ \ell, m}.
\end{equation*}
In particular, 
$\mathcal {N}_m^2$
is $\mathcal {N}_m$.  

Clause (4) of the theorem 
is immediate from the special case $j=2$  of the following claim,
since $\mathbf N^{2, m}$ belongs to 
$\mathcal {N}_m^2 =\mathcal {N}_m$,
and since each $\mathcal {N}_m^j $ 
is a subvariety of $ \mathcal {N}_m$.

\begin{claim} \labbel{claim}  
For every $j$ such that  $2 \leq j \leq \ell$,
there are an algebra $\mathbf A_3^j \in \mathcal {N}_m^j$
and a subalgebra $\mathbf F^j$ of 
$ \mathbf A_3^j \times \mathbf N^{2, m} $ 
such that the congruence identity 
\begin{equation}\labbel{ciddd}
\alpha (\beta \circ ( \alpha \gamma \circ \alpha \beta 
\circ {\stackrel{q-2}{\dots}} \circ \alpha \beta ^\bullet   )
\circ   \gamma ^ \bullet) 
\subseteq
(\alpha ( \gamma \circ \beta \circ {\stackrel{q}{\dots}}
 \circ  \beta ^ \bullet  )) ^{m-2j+2}  
   \end{equation}    
fails in $\mathbf F^j$.
\end{claim}    

In order to prove the claim
we need to establish some  notation.

\begin{notation} \labbel{nota}    
Let  $\mathbf C_{q+1}$
be the chain with $q+1$
elements $\{ 0, 1, \dots, q \}$ 
and the standard lattice operations. 
Let $  \mathbf N _{q+1}^{j,m} $ denote
$   \mathbf C_{q+1} ^{{\rm nu},j,m}$,
that is,  recalling Definition \ref{vars}, 
 $  \mathbf N _{q+1}^{j,m} $ is the term-reduct of 
$\mathbf C_{q+1}$ with the only operation given by 
the term $u_{j,m}$ from Definition \ref{latterms}.  
In particular,
$  \mathbf N ^{j,m} $ is $  \mathbf N_2 ^{j,m} $.

For every $q \geq 2$, let
$\beta^*_{q+1}$
be the congruence on 
$\mathbf C_{q+1}$
determined by the partition
$\{ \{ q,q-1\}, \{ q-2, q-3\}, \dots  \}$,
where $\{ 0 \}$ is a block of 
$\beta^*_{q+1}$ if $q$ is even. Let  
 $\gamma^*_{q+1}$  
be the congruence on 
$\mathbf C_{q+1}$
determined by the partition
$\{ \{ q\}, \{ q-1, q-2\}, \{ q-3, q-4\}, \dots  \}$,
where $\{ 0 \}$ is a block of 
$ \gamma ^*_{q+1}$ if $q$ is odd.
Notice that $\beta^*_{q+1}$
and $\gamma^*_{q+1}$ are congruences 
on every term-reduct of   
$\mathbf C_{q+1}$.

 If $\mathbf A$ is an algebra, we
let $0_{\mathbf A}$  
 denote
the smallest congruence on $\mathbf A$.
Similarly, $1_{\mathbf A}$
denotes
the largest congruence on $\mathbf A$.
When  there is 
no risk of ambiguity we shall omit  subscripts.
\end{notation}

The claim is proved 
in three steps
by induction on decreasing $\ell$.
During the inductive proof of the claim
we shall need some further properties of the  constructions
witnessing the claim itself.
Recall that   
$\mathbf N^{2, m}$
is a reduct of the two-elements lattice with 
base set $\{ 0,1\}$. 
We shall need the additional properties
stated in the following subclaim.

\begin{subclaim}    
(*) For every $j$, the failure of 
identity \eqref{ciddd} in
$\mathbf F^j$ can be witnessed by elements
of the form
$(a^j,1) $ and $(d^j,1)$.
By the above statement we mean that
we can choose $a^j$ and $d ^j$ in $\mathbf A_3^j$ 
and congruences $\alpha$, $\beta$ and $\gamma$ of $\mathbf F^j$
in such a way that
$(a^j,1) ,(d^j,1) \in F^j$ and
  the pair $((a^j,1) ,(d^j,1))$
belongs to the left-hand side of 
\eqref{ciddd}, but not to the right-hand side
of \eqref{ciddd}.

(**) We shall also require that the assumption 
in (*) above
that the pair 
$((a^j,1) ,(d^j,1))$ belongs to 
$\beta \circ ( \alpha \gamma \circ \alpha \beta 
\circ {\stackrel{q-2}{\dots}} \circ \alpha \beta ^\bullet)
\circ \gamma ^\bullet$ 
can be witnessed by elements of the form
$(c_i^j,0)$, namely, that there are elements 
$c_1^j, \dots, c_{q-1}^j$ in $A_3^j$ such that 
$(a^j,1) \mathrel { \beta } (c^j_1,0)
\mathrel { \alpha \gamma } (c^j_2,0)
\mathrel { \alpha \beta } (c_3^j,0)
\mathrel { \alpha \gamma}  \dots 
\mathrel { \alpha \beta ^\bullet } (c^j_{q-1},0)
\mathrel {  \gamma  ^\bullet} (d^j,1)$
and $(c^j_1,0), (c^j_2,0), \dotsc \in F^j$. 

(***) Finally, we shall prove that we can make 
 \eqref{ciddd} fail by taking 
 $\alpha$  
to be the congruence induced on $\mathbf F^j$
by the congruence 
$ 1   \times  0 $
on  
$ \mathbf A_3^j \times \mathbf N^{2, m} $.
Actually, we shall only need that 
the second component is 
$0$, but the proof shall give the  
additional result on the first component.
\end{subclaim}

We now proceed with the proof of the claim,
at the same time checking that we can handle the proof 
in such a way that (*) - (***) in the subclaim
are  verified.

 \smallskip 

\emph{First step.} Consider the case 
when $m$ is odd and $j= \ell$, thus 
$m-2j+2 = 1$.
In this case the claim is almost obvious since 
if the exponent on the right is $1$, then 
 identity 
\eqref{ciddd} implies congruence $q$-permutability
(just take $\alpha=1$, the largest congruence).
Lattices are not $q$-permutable,
hence, a fortiori,  the term-reduct
$\mathcal {N}_m^\ell$ 
is not $q$-permutable.
It is then enough to take some witness
$\mathbf A_3^\ell \in \mathcal {N}_m^\ell$
of the failure of $q$-permutability 
and take $\mathbf F^\ell = \mathbf A_3^\ell 
\times \mathbf N^{2, m} $.

 In detail,
recall Notation \ref{nota}, 
take
 $\mathbf A_3^\ell
=    \mathbf N_{q+1} ^{\ell,m} $,
$\mathbf F^\ell = \mathbf A_3^\ell 
\times \mathbf N^{2, m} $
and consider
the elements
$(q,1), (q-1,0),  (q-2,0), \dots, (1,0), (0,1)$ 
and the congruences
$ \beta = \beta^*_{q+1}  \times 1 $,   
$ \gamma = \gamma ^*_{q+1} \times 1 $
and 
$\alpha= 1 \times 0$ on $\mathbf F^\ell $
in order to get the failure of \eqref{ciddd}
and at the same time to have (*) - (***) 
satisfied.
Indeed, the pair 
$((q,1),  (0,1))$ 
does belong to 
$\alpha (\beta \circ ( \alpha \gamma \circ \alpha \beta 
\circ {\stackrel{q-2}{\dots}} \circ \alpha \beta ^\bullet   )
\circ   \gamma ^ \bullet) $,
as witnessed by the other elements in the above sequence.
On the other hand, suppose by contradiction that 
$((q,1),  (0,1))$  belongs to
$\alpha ( \gamma \circ \beta \circ {\stackrel{q}{\dots}}
 \circ  \beta ^ \bullet  ) $.
Then $(q,1) \mathrel { \gamma   } (e_1, f_1) \mathrel { \beta  }
(e_2, f_2) \mathrel { \gamma  } \dots   \mathrel { \gamma ^\bullet  }
(e_{q-1}, f_{q-1}) 
 \mathrel { \beta ^\bullet  } (0,1)$, for certain elements  
 $(e_1, f_1), (e_2, f_2), \dotsc  \in 
A_3^\ell \times N^{2, m}$. 
By $ \gamma $-equivalence, 
$e_1=q$; then, by $ \beta $-equivalence, 
$e_2 \geq q-1$ and, again by $ \gamma $-equivalence,
 $e_3 \geq q-2$. Going on, 
$e_{q-1} \geq 2$, hence
$(e_{q-1}, f_{q-1})$ is not   
$\beta ^\bullet $-equivalent to $(0,1)$.

Notice that $\mathbf N_{q+1} ^{\ell,m}  $
belongs to the variety generated by 
 $ \mathbf N ^{\ell,m} = \mathbf N_2 ^{\ell,m}$,
since  $\mathbf C_{q+1}$ belongs to the variety 
generated by $\mathbf C_{2}$.
In particular, $\mathbf A_3^\ell =    \mathbf N_{q+1} ^{\ell,m}$  belongs
to $\mathcal {N}_m^{\ell}$. 

\smallskip

\emph{Second step.} Next,  we consider the case $j= \ell$ and 
 $m$ even in the claim. In this case 
$m-2j+2 = 2$.
Apply Lemma \ref{lemnu}
taking  
$\mathbf A_1=\mathbf A_2=
\mathbf N_{q+1} ^{\ell,m} $ and
  $\mathbf A_3$ a one-element algebra
with an $m$-ary operation, say, $A_3= \{  a\} $, with $a=d$.
 Finally, let
$\mathbf A_4 = \mathbf N ^{2,m} =  \mathbf N_2 ^{2,m}$
and let $\mathbf F$ be  the whole of 
$ \mathbf A_3 \times  \mathbf A_4 $.

Take $h=k=\ell$ in Lemma \ref{lemnu}. 
By Observation \ref{obs}, the element $0$ is $ \ell $-absorbing in   
$\mathbf A_1 $ and in $ \mathbf A_2$
and $0_4$ is $2$-absorbing in  $ \mathbf A_4$.
The operation of  $\mathbf A_3$
is trivially an
$ \ell $-majority term.
By Lemma \ref{lemnu}
we get a subalgebra 
$\mathbf B$ of  
$\mathbf A_1 \times \mathbf A_2 \times\mathbf A_3\times \mathbf A_4$.

The  proof that \eqref{ciddd} fails in $\mathbf B$
for $j= \ell$
presents no significant difference with respect to 
\cite{B,day}. We recall the details.
Consider the following elements of $B$.
\begin{equation}\labbel{form}
\begin{aligned}
c_0 & =(q,0,a,1),  \qquad && c_q=  (0,q,a,1),
\qquad  \text{ and } 
\\
 c_i&= (q-i,i, a,0),  \qquad  && \text{for $i=1, \dots,q-1$.} 
\end{aligned}
  \end{equation}     
The above elements are indeed in $B$,
since   $c_0$ has type
I$^ \sigma$, $ c_q$  has type
III$^ \sigma$ (since $a=d$) and
the remaining 
$  c_i$'s
have type
IV$^ \sigma$.
Recall that in this special case we have taken 
$\mathbf F $ equal to $  \mathbf A_3\times \mathbf A_4$,
hence the above elements automatically  belong to
$A_1 \times A_2 \times F$. 

If $q$ is even, let $  \beta$ and $  \gamma$ be, respectively, 
the congruences on $\mathbf B$ induced by
$\beta^*_{q+1} \times \gamma^*_{q+1}  \times  1 \times 1 $ and   
$ \gamma ^*_{q+1} \times \beta ^*_{q+1}  \times   1 \times 1 $.
If $q$ is odd, let $  \beta$ and $  \gamma$ be, respectively, 
the congruences on $\mathbf B$ induced by
$\beta^*_{q+1} \times \beta ^*_{q+1}  \times     1 \times 1$ and   
$ \gamma ^*_{q+1} \times \gamma  ^*_{q+1}    \times   1 \times 1$.
Both in case $q$ even and  $q$ odd, let 
 $  \alpha $ be
the congruence induced by
$1 \times 1   \times   1 \times 0 $.

We have $ c_0 \mathrel {    \alpha }  c_q $
and $ c_0 \mathrel {    \beta  }    c_1 
 \mathrel {    \alpha    \gamma  }    c_2
 \mathrel {    \alpha    \beta   }    c_3 \dots$,
hence  
$ (c_0,   c_q ) 
\in 
\alpha (\beta \circ ( \alpha \gamma \circ \alpha \beta 
\circ {\stackrel{q-2}{\dots}} \circ \alpha \beta ^\bullet   )
\circ   \gamma ^ \bullet) $.
We shall show that 
$ (c_0,   c_q ) 
\notin 
\alpha ( \gamma \circ \beta \circ {\stackrel{q}{\dots}}
 \circ  \beta ^ \bullet  ) \circ
\alpha ( \gamma \circ \beta \circ {\stackrel{q}{\dots}}
 \circ  \beta ^ \bullet  )$ in $\mathbf B$.
Towards a contradiction, suppose the contrary.
Then there is some element $f \in B$ 
such that 
$ (c_0,   f) 
\in 
\alpha ( \gamma \circ \beta \circ {\stackrel{q}{\dots}}
 \circ  \beta ^ \bullet  )$
and 
$ (f,c_q) 
\in 
\alpha ( \gamma \circ \beta \circ {\stackrel{q}{\dots}}
 \circ  \beta ^ \bullet  )$.
Thus
$c_0 \mathrel { \alpha } f $
and there are elements 
$f_0=c_0, f_1, \dots, f_q=f$
such that 
$f_0 \mathrel \gamma \allowbreak f_1 \mathrel \beta \allowbreak   f_2 \dots$   
Recall that 
$f_0=c_0= (q,0,a,1)$.
By $ \gamma $-equivalence of
$f_0$ and  $f_1$, 
the first component of $f_1$
 is $q$. By $\beta$-equivalence 
of $f_1$ and  $f_2$,
the first component of $f_2$
is  $ \geq q-1$.
Going on, the first component
of $f_q=f$ is  $\geq 1$,
in particular, it is not $0$.
Thus $f$ has neither type II$^ \sigma$ nor III$^ \sigma$.
Moreover,
$f$ has not type IV$^ \sigma$, either,
since its fourth component is $1$, 
by $\alpha$-equivalence
of $f_0$ and $f$.
Since $f \in B$, then 
$f$ has necessarily type I$^ \sigma$, thus its second component is $0$.
However, by performing a symmetric argument,
using the assumption
$ (f,c_q) 
\in 
\alpha ( \gamma \circ \beta \circ {\stackrel{q}{\dots}}
 \circ  \beta ^ \bullet  )$,
we get that the second component of $f$
is $\geq 1$, a contradiction.  
We have showed that \eqref{ciddd} fails in $\mathbf B$
for $j= \ell$.

We are almost done. It is now enough to 
 declare who $\mathbf A_3^\ell$ and 
$\mathbf F^\ell$ actually are. 
Take $\mathbf A_3^\ell$ 
to be $\mathbf A_1 \times \mathbf A_2 \times \mathbf A_3$.
As in the first step, 
$\mathbf A_1=\mathbf A_2=
\mathbf N_{q+1} ^{\ell,m} $  belong
to the variety $\mathcal {N}_m^{\ell}$,
hence $\mathbf A_3^\ell$ belongs to 
$\mathcal {N}_m^{\ell}$, too.
Finally, let $\mathbf F^\ell = \mathbf B$,
thus  \eqref{ciddd} fails in 
$\mathbf F^\ell$.
The additional conditions (*) - (***)
are verified by  construction.
Indeed, here we take
$a^j = (q,0,a)$
and $d^j=(0,q,a)$,
thus  
$c_0  =(a^j,1)$
and $c_q  =(d^j,1)$,
modulo a standard identification
of nested components.
Similarly for the other $c_i$'s.

Let us observe that in the present step we could have 
worked  with just three coordinates.
However, it is easier to add a dummy third coordinate, 
rather than state and prove also a three-coordinate 
(and essentially less general)
version
of Lemma \ref{lemnu}. 
The full four-coordinate version of 
Lemma \ref{lemnu} will be necessary in the proof of the
next step.

 \smallskip 

\emph{Third step.}
Finally,  we suppose that we have proved the claim 
and the subclaim for some 
$j$ with 
$2 <  j \leq \ell$
and we shall prove the claim 
for $j-1$. 
Since we have proved the claim when 
$j=\ell$, 
 an easy finite induction establishes the claim
and the subclaim for all $j$'s,
hence the theorem.

The proof of the third step is not
really different from the proof of the second step.
However, here we shall use a nontrivial  $\mathbf A_3$
which is given by the inductive hypothesis.
Taking into account a nontrivial $\mathbf A_3$ 
involves a bit of further details and, as we mentioned,
the full power of Lemma \ref{lemnu}
will be necessary. 

So let $\mathbf A_3^j $
and $\mathbf F^j \subseteq \mathbf A_3^j \times \mathbf N^{2,m}  $
 be given by the case $j$ of the claim
and let the failure of 
\eqref{ciddd} in 
$\mathbf F^j$ be witnessed by congruences
$ \tilde \alpha$, $ \tilde \beta$ and $\tilde \gamma$.
We can inductively assume that
properties (*) - (***) hold,
so, by (*),  let
$((a^j, 1) ,(d^j,1))$
belong to the left-hand side of 
\eqref{ciddd}, but not to the right-hand side,
for certain $a^j,d^j \in \mathbf A_3^j $ and  where
$\alpha$, $\beta$ and $\gamma$ in \eqref{ciddd}
are replaced by
$ \tilde \alpha$, $ \tilde \beta$ and $\tilde \gamma$.

Apply Lemma \ref{lemnu}
taking $h=j-1$, $k=m{-}j{+}1$, 
$\mathbf A_1=\mathbf A_2=
\mathbf N_{q+1} ^{j-1,m} $, 
$\mathbf A_3 =\mathbf A_3^j $,
$\mathbf A_4 = \mathbf N ^{2,m} $,
$a=a^j$, $b=b^j$    
and $\mathbf F = \mathbf F^j$. 
Again, the algebra $ \mathbf N_{q+1} ^{j-1,m} 
=\mathbf C_{q+1} ^{{\rm nu},j-1,m} $
belongs to the variety generated by 
 $ \mathbf N ^{j-1,m} =  \mathbf C_{2} ^{{\rm nu},j-1,m}$,
since  $\mathbf C_{q+1}$ belongs to the variety 
generated by $\mathbf C_{2}$.
In particular, $\mathbf A_1=\mathbf A_2=
\mathbf N_{q+1} ^{j-1,m} $ belong
to $\mathcal {N}_m^{j-1}$. 

By Observation \ref{obs}, $0$ is $j-1$-absorbing in   
$\mathbf A_1 $ and in $ \mathbf A_2$.
Moreover, the operation of  $\mathbf A_3 $
is an
$ m{-}j{+}1$-majority term,
since  $\mathbf A_3 =\mathbf A_3^j \in \mathcal {N}_m^{j}$ and
each operation on the  generators of
$\mathcal {N}_m^{j}$ is an  $ m{-}j{+}1$-majority term,
again by Observation \ref{obs}.

Hence we can apply Lemma \ref{lemnu}
with $h=j-1$   and $k=  m{-}j{+}1$
(notice that $h \leq k$, since $j \leq \ell$),
getting a subalgebra 
$\mathbf B$ of  $\mathbf A_1 \times \mathbf A_2 \times \mathbf F$,
which is itself  a subalgebra of 
$\mathbf A_1 \times \mathbf A_2 \times\mathbf A_3\times \mathbf A_4$.

Recall the definitions of 
 $\beta^*_{q+1} $ and $ \gamma^*_{q+1}$
from Notation \ref{nota}.
If $q$ is even, let $\beta$ and $ \gamma$ be, respectively, 
the congruences on $\mathbf B$ induced by
$\beta^*_{q+1} \times \gamma^*_{q+1}  \times \tilde   \beta $ and   
$ \gamma ^*_{q+1} \times \beta ^*_{q+1}  \times  \tilde   \gamma $.
If $q$ is odd, let $ \beta$ and $ \gamma$ be, respectively, 
the congruences on $\mathbf B$ induced by
$\beta^*_{q+1} \times \beta ^*_{q+1}  \times \tilde   \beta$ and   
$ \gamma ^*_{q+1} \times \gamma  ^*_{q+1}  \times \tilde    \gamma $.
In both cases, let 
 $\alpha $ be
the congruence induced by
$1 \times 1   \times \tilde  \alpha  $.
By (**), 
 there are elements 
$   c_1^j, \dots,    c^j_{q-1}$ in $A_3^j$ such that 
$(a^j,1) \mathrel { \tilde \beta } ( c_1^j, 0)
\mathrel {\tilde \alpha \tilde \gamma } (  c_2^j, 0) 
\mathrel {\tilde \alpha \tilde \beta  } (  c_3^j, 0) \dots $  
Consider the following elements of $B$.
\begin{align*}
 c_0& =(q,0,a^j,1),   & c_{q} = (0,q,d^j, 1),
\qquad  \text{ and } 
\\
 c_i&= (q-i,i,   c^j_i, 0),   &     \text{for $i=1, \dots,q-1$.} \qquad
\end{align*} 
  Notice that $c_0$ belongs to 
$A_1 \times A_2 \times F$,
since $(a^j,1) \in F$, by (*).
Moreover,  $ c_0$ has type
I$^ \sigma$, hence $c_0 $ is indeed in $  B$.
Recall that we are taking $a=a^j$. 
Similarly, 
$ c_q$ belongs to 
$A_1 \times A_2 \times F$ and has type
III$^ \sigma$, hence $ c_q \in B$.
The remaining 
$ c_i$'s belong to 
$A_1 \times A_2 \times F$,
since $( c_i^j, 0) \in F$, by (**).
Moreover, each $ c_i$ 
has type
IV$^ \sigma$, hence $c_i \in B$.

One easily checks that 
$ (c_0, c_q)\in 
 \alpha ( \beta \circ (  \alpha  \gamma \circ  \alpha  \beta 
\circ {\stackrel{q-2}{\dots}} \circ  \alpha  \beta ^\bullet   )
\circ    \gamma ^ \bullet)$.
We shall show that
$ (c_0, c_q)\notin
 ( \alpha (  \gamma \circ  \beta \circ {\stackrel{q}{\dots}}
 \circ   \beta ^ \bullet  )) ^{m-2j+4} $,
thus identity \eqref{ciddd}  
fails in $\mathbf B$ for $j-1$. 
Suppose the contrary.
  Then there are elements 
$f,g \in B$ such that 
$ (  c_0,  f)\in
  \alpha (  \gamma \circ  \beta \circ {\stackrel{q}{\dots}}
 \circ   \beta ^ \bullet  ) $,
 $ (f, g)\in
 ( \alpha (  \gamma \circ  \beta \circ {\stackrel{q}{\dots}}
 \circ   \beta ^ \bullet  )) ^{m-2j+2} $
and
$ (g,   c_q)\in
  \alpha (  \gamma \circ  \beta \circ {\stackrel{q}{\dots}}
 \circ   \beta ^ \bullet  ) $.
Notice that 
$m-2j+2 \geq 1$, since $j \leq \ell$. 
From the first relation
we get that  $  c_0 \mathrel {    \alpha }  f $ and 
that there are
elements 
$    c_0=f_0, f_1, \dots, f_q=  c_q$ in $B$ 
such that $f_0 \mathrel {     \gamma  } f_1 \mathrel {     \beta   } f_2
\mathrel {    \gamma   } f_3 \dots  $
Since 
$  c_0= f_0 =(q,0,a^j,1)$,
then, by $    \gamma$-equivalence, the first component of $f_1$ is $q$.      
By $    \beta $-equivalence, the first component of $f_2$ is   $ \geq q-1$.
Going on,   the first component of $f_q=f$ is $ \geq 1$,
thus $f$ has neither type II$^ \sigma$ nor type III$^ \sigma$.  
Since $  c_0 \mathrel {    \alpha }  f $,
then, by (***), the fourth component
of $f$ is $1$, hence $f$ has not type IV$^ \sigma$, either. Since $f \in B$,
then $f$ has type I$^ \sigma$, thus the third component of $f$ is $a^j$.     
Symmetrically, the fourth component 
of $g$ is $1$,  $g$ has type III$^ \sigma$, hence the third component of $g$
is $d^j$. From $ (f, g)\in
 ( \alpha (  \gamma \circ  \beta \circ {\stackrel{q}{\dots}}
 \circ   \beta ^ \bullet  )) ^{m-2j+2} $,
restricting to the third  component
of $\mathbf A_1 \times \mathbf A_2 \times\mathbf F$,
we get 
$ ((a,^j,1),(d^j,1))\in
 ( \tilde  \alpha ( \tilde   \gamma \circ \tilde   \beta \circ {\stackrel{q}{\dots}}
 \circ  \tilde   \beta ^ \bullet  )) ^{m-2j+2} $,
contradicting our assumption that 
the pair $((a^j,1),(d^j,1))$
witnesses the failure of 
\eqref{ciddd} for 
$ \tilde   \alpha$, $ \tilde  \beta$ and
 $ \tilde  \gamma$. 

We have showed 
that $ (  c_0,   c_q)\in 
 \alpha ( \beta \circ (  \alpha  \gamma \circ  \alpha  \beta 
\circ {\stackrel{q-2}{\dots}} \circ  \alpha  \beta ^\bullet   )
\circ    \gamma ^ \bullet)$ and 
$ (  c_0,   c_q)\notin
 ( \alpha (  \gamma \circ  \beta \circ {\stackrel{q}{\dots}}
 \circ   \beta ^ \bullet  )) ^{m-2j+4} $.
Now it is enough to
take 
$\mathbf A_3^{j-1} = \mathbf A_1 \times \mathbf A_2 \times \mathbf A_3^{j}$
and $\mathbf F^{j-1}= \mathbf B \subseteq \mathbf A_3^{j-1} 
\times \mathbf N^{2,m} $ (the inclusion
is considered modulo the usual identifications), 
to get that  \eqref{ciddd} fails in $\mathbf F^{j-1}$ for $j-1$.
Notice that $\mathbf A_3^{j}$ belongs to $\mathcal {N}_m^{j}$,
by the inductive assumption.
As we mentioned, 
$ \mathbf A_1 $ and $  \mathbf A_2 $
belong to $\mathcal {N}_m^{j-1}$,
hence
$\mathbf A_3^{j-1}$
belongs to $\mathcal {N}_m^{j-1}$, too,
since 
$\mathcal {N}_m^{j} \subseteq \mathcal {N}_m^{j-1}$.
As in the second step, (*) - (***) are verified by construction.
\end{proof}

Theorem \ref{mmm} is optimal.
It is immediate from Theorem \ref{mnun}
that  clauses (1) and (3) in \ref{mmm}
are optimal. Clause (4), too, is the best possible result,
as shown in the next remark.

\begin{remark} \labbel{mult} 
In \cite[Proposition 5.1]{B} 
we have showed that if $m \geq 3$
and  some variety $\mathcal {V}$ 
has an $m$-ary near-unanimity term, then,
for every $q \geq 2$,
$\mathcal {V}$ satisfies
\begin{align} \labbel{bah}      
\alpha  (\beta \circ  \gamma \circ {\stackrel{q}{\dots}} \circ  \gamma  )
&\subseteq 
\alpha  \beta  \circ \alpha \gamma \circ 
 {\stackrel{(m-2)q}{\dots}} \circ  \alpha  \gamma, 
&&\text{if $q$ is even,}
\\ 
\labbel{bahbah}
\alpha  (\beta  \circ \gamma \circ {\stackrel{q}{\dots}} \circ  \beta   )
 &\subseteq 
\alpha \beta   \circ  \alpha \gamma 
\circ {\stackrel{1+(m-2)(q-1)}{\dots\dots}} \circ  \alpha  \beta ,
&& \text{if $q$ is  odd.} 
\end{align}    

Clause (4) 
in Theorem \ref{mmm} shows that the above result is best possible.
Indeed, $\mathcal {N}_m$ has an $m$-ary near-unanimity term.
If $q$ is even
and, by contradiction, \eqref{bah}
can be improved by considering $(m-2)q-1$
factors on the right-hand side, then  
\begin{align*} 
\alpha (\beta \circ ( \alpha \gamma \circ \alpha \beta 
\circ {\stackrel{q-2}{\dots}} \circ \alpha \beta  )
\circ   \gamma ) 
&\subseteq
\alpha  (\beta \circ  \gamma \circ \beta \circ
 {\stackrel{q}{\dots}} \circ \beta  \circ  \gamma  ) 
\\
\subseteq  \alpha  \beta  \circ \alpha \gamma \circ 
 {\stackrel{(m-2)q-1}{\dots}} \circ  \alpha  \beta  
&\subseteq  
\alpha \gamma \circ \alpha  \beta  \circ \alpha \gamma \circ 
 {\stackrel{(m-2)q}{\dots}} \circ  \alpha  \beta
\\
& \subseteq 
(\alpha ( \gamma \circ \beta \circ {\stackrel{q}{\dots}}
 \circ  \beta   )) ^{m-2},  
\end{align*} 
contradicting clause (4) in Theorem \ref{mmm}.
 
The case $q$ odd is similar, using  
Lemma \ref{cilem} and identity \eqref{ciodd}.  
In fact, the arguments show that
if $m \geq 3$, then   $\mathcal {N}_m$  is  a variety
 with a symmetric $m$-ary near-unanimity term
for which  the following identities fail.
\begin{align*}      
\alpha  (\beta \circ  \gamma \circ {\stackrel{q}{\dots}} \circ  \gamma  )
&\subseteq 
\alpha  \gamma   \circ \alpha \beta  \circ 
 {\stackrel{(m-2)q}{\dots}} \circ  \alpha  \beta , 
&&\text{if $q$ is even,}
\\ 
\alpha  (\beta  \circ \gamma \circ {\stackrel{q}{\dots}} \circ  \beta   )
 &\subseteq 
\alpha \gamma    \circ  \alpha \beta  
\circ {\stackrel{1+(m-2)(q-1)}{\dots\dots}} \circ  \alpha  \gamma  ,
&& \text{if $q$ is  odd.} 
\end{align*}
Notice that here $ \alpha \beta $ and
$ \alpha \gamma $ are exchanged on the right-hand side,
in comparison with \eqref{bah}.
\end{remark}

\section{An explicit example} \labbel{expl}

Following the proof of Theorem \ref{mmm}
we shall  present explicit examples of algebras 
in $\mathcal {N}_m$ for which identity \eqref{cid}
fails.  Then we hint to the details for a full counterexample
to \eqref{cidd}. 
In particular, such counterexamples show 
in a more direct way that 
Theorem \ref{mnun} cannot be improved. 
  
Recall that 
$\mathbf C_{q+1}$ is the chain with 
$q+1$ elements $\{ 0, 1, \dots, q\}$,
considered as a lattice.
Moreover,  
$\mathbf N_{q+1}^{ j,m}$ is 
the term-reduct of 
$\mathbf C_{q+1}$ endowed with the  $m$-ary operation 
induced by the lattice term 
$u _{j,m}(x_1, \dots, x_m)= \prod _{|J|=j} \sum _{i \in J} x_i$,
 where $J$ varies on  subsets of 
$\{ 1, \dots, m \}$.
We have set  $\mathbf N^{ j,m} = \mathbf N_{2}^{ j,m}$.
Compare Definitions \ref{latterms}, \ref{vars}   and Notation \ref{nota}. 

Fix some $q \geq 2$, $m \geq 3$ and 
let $\ell= \frac{m}{2} $ if 
$m$ is even and
 $\ell= \frac{m+1}{2} $ if 
$m$ is odd.
 Consider the following product $\mathbf P = \mathbf P (m,q)$  in the cases, respectively, $m$
even and $m$ odd.  
\begin{equation}\labbel{pro} \tag{P}
\begin{gathered}     
(\mathbf N_{q+1} ^{2,m} 
\timee
\mathbf N_{q+1} ^{2,m} )
\timee
(\mathbf N_{q+1} ^{3,m} 
\timee
\mathbf N_{q+1} ^{3,m} )
\timee
{\hspace {0.4pt} \dots \hspace {0.4pt}}
\timee
(\mathbf N_{q+1} ^{\ell,m} 
\timee
\mathbf N_{q+1} ^{\ell,m} )
\timee
\mathbf N ^{ 2,m},
\\ 
(\mathbf N_{q+1} ^{2,m} 
\timee
\mathbf N_{q+1} ^{2,m} )
\timee
(\mathbf N_{q+1} ^{3,m} 
\timee
\mathbf N_{q+1} ^{3,m} )
\timee
{ \dots }
\timee
(\mathbf N_{q+1} ^{\ell{-}1,m} 
\timee
\mathbf N_{q+1} ^{\ell{-}1,m} )
\timee
\mathbf N_{q+1} ^{\ell,m}
\timee
\mathbf N ^{ 2,m},
\end{gathered}
   \end{equation}     
where the grouping of the factors is only for notational convenience.
In any case, $ \mathbf P (m,q)$ has $m-1$ factors;
for example,  $ \mathbf P (3,q) = \mathbf N_{q+1} ^{2,3}
\times
\mathbf N ^{ 2,3}$,
$ \mathbf P (4,q) = (\mathbf N_{q+1} ^{2,4} 
\times
\mathbf N_{q+1} ^{2,4} )
\times
\mathbf N ^{ 2,4}$
and 
$ \mathbf P (5,q) = 
(\mathbf N_{q+1} ^{2,5} 
\times
\mathbf N_{q+1} ^{2,5}) \times
\mathbf N_{q+1} ^{3,5}
\times
\mathbf N ^{ 2,5}$.
By Observation \ref{obs}, the operation of $\mathbf P$  
is an $m$-ary near-unanimity term. 

A member $p$ of $P$ is \emph{good}
if either (a) its last component is $0$, or (b) its last component is
$1$ and, disregarding the last component,  

(b1) $p$  begins with a (possibly empty, possibly covering all pairs)
sequence of null  pairs $(0,0)$,

(b2) the first (if any) pair of $p$ which is not null   
has either the form $( \hyphe, 0)$, or the form  
$( 0, \hyphe)$, and

(b3) all the subsequent pairs, if any, have, correspondingly,
the form $(q,0)$ or  $(0,q)$.

If $m$ is odd, we follow the same rules, considering
the penultimate component as a ``half pair'' and applying the above rules
only to the first component of the pair.   
Typical good elements are
given by the sequences 
\begin{equation*}
\begin{tabular}{cccccccccc}
$(\hyphe,\hyphe)$ &$(\hyphe,\hyphe)$ & $\dots$ & 
 $(\hyphe,\hyphe)$ & $( \hyphe, \hyphe)$ & $ (\hyphe, \hyphe)$ & $\dots$ & 
 $(\hyphe,\hyphe) $ & $\hyphe$ & $ 0 $
\\
$( \hyphe, 0)$ &$(q,0)$ & $\dots$ & 
 $(q,0)$ & $( q, 0)$ & $ (q, 0)$ & $\dots$ & 
 $(q,0) $ & $q $ &$ 1 $
\\
$(0,0)$ &$(0,0)$ & $\dots$ & 
 $(0,0)$ & $( \hyphe, 0)$ & $ (q, 0)$ & $\dots$ & 
 $(q,0) $ & $ q$ &$ 1 $
\\
$(0,0)$ &$(0,0)$ & $\dots$ & 
 $(0,0)$ & $( 0, 0)$ & $ (0, 0)$ & $\dots$ & 
 $(\hyphe,0) $ & $ q$ &$ 1 $
\\
$(0,0)$ &$(0,0)$ & $\dots$ & 
 $(0,0)$ & $( 0, 0)$ & $ (0, 0)$ & $\dots$ & 
 $(0,0) $ & $ \hyphe$ &$ 1 $
\\
$(0,\hyphe)$ &$(0,q)$ & $\dots$ & 
 $(0,q)$ & $(0, q)$ & $ (0,q)$ & $\dots$ & 
 $(0,q) $ & $0$ &$ 1 $
\\
$(0,0)$ &$(0,0)$ & $\dots$ & 
 $(0,0)$ & $(0, \hyphe)$ & $ (0,q)$ & $\dots$ & 
 $(0,q) $ &  $0$ &$ 1 $
\\
$(0,0)$ &$(0,0)$ & $\dots$ & 
 $(0,0)$ & $( 0, 0)$ & $ (0, 0)$ & $\dots$ & 
 $(0,\hyphe) $ & $ 0$ &$ 1 $
\end{tabular}
  \end{equation*}    
in the case $m$ odd,
while we get typical elements in the case $m$
even simply discarding the penultimate column
in the above table. 
 
The set of good elements of $P$ 
is the universe for a subalgebra $\mathbf B = \mathbf B (m,q)$ of $\mathbf P$.
This can be checked directly using arguments similar to 
those used in the proof of Lemma \ref{lemnu}.
Roughly, suppose that  $b_1, \dots, b_m \in B$
and $b=u(b_1, \dots, b_m)$. If at least two $b_i$'s
have $0$ as the last component, then this applies to 
$b$, as well, hence $b \in B$, by (a). 
Otherwise, there are enough $0$'s in the components 
of the $b_i$'s in order 
to make $0$ at least one element of each pair of $b$,
using Observation \ref{obs}.
Then the rules describing the elements of $B$,
together  with Observation \ref{obs} again,
show that $b$  has  a sufficient number of $0$'s and $q$'s
in  the appropriate places. 
Alternatively, in order to show
that $B$  is the universe for a subalgebra  of $\mathbf P$,
 work out the proof of Theorem \ref{mmm},
going in the backward direction.    

Now suppose, say, that  $q$ is even
and consider the congruences  $\beta$, 
$\gamma$ and $\alpha$  on  $\mathbf B$
induced, respectively, by the congruences 
\begin{equation}\labbel{gath}
   \begin{gathered}
\beta ^*  = 
(\beta^*_{q+1} 
\times
\gamma ^*_{q+1} )
 \times
(\beta^*_{q+1} 
\times
\gamma ^*_{q+1} )
\times
\dots
\times
(\beta^*_{q+1} 
\times
\gamma ^*_{q+1} )
\times
\beta^*_{q+1}
\times
1,
\\
 \gamma ^*   = 
(\gamma^*_{q+1} 
\times
\beta ^*_{q+1} )
\times
(\gamma^*_{q+1} 
\times
\beta ^*_{q+1} )
\times
\dots
\times
(\gamma^*_{q+1} 
\times
\beta ^*_{q+1} )
\times
\gamma^*_{q+1}
\times
1,
\\
\alpha ^*   = 
(1  
\times
1  )
\times
(1  
\times
1  )
\times
\dots
\times
(1  
\times
1  )
\times
1 
\times
0,
\end{gathered}
  \end{equation}     
where, as usual by now, in each line the 
penultimate congruence appears only if $m$ is odd.
Recall that
$\beta^*_{q+1}$ and 
 $\gamma^*_{q+1}$  
are the congruences 
determined, respectively, by the partitions
$\{ \{ q,q-1\}, \{ q-2, q-3\}, \dots  \}$ and
$\{ \{ q\}, \{ q-1, q-2\}, \{ q-3, q-4\}, \dots  \}$.

For simplicity, let $q=2$. Indeed, this is the case 
showing that Theorem \ref{mnun}(1) cannot be improved.
Consider the following elements of $P$  
\begin{equation*}
\begin{tabular}{lccccccccccc}
$a $ & $ =$ & $( 2, 0)$ &$(2,0)$ & $\dots$ & 
 $(2,0)$ & $( 2, 0)$ & $ (2, 0)$ & $\dots$ & 
 $(2,0) $ & $2 $ &$ 1 $
\\
$c $ & $ =$ & $(1,1)$ &$(1,1)$ & $\dots$ & 
 $(1,1)$ & $( 1, 1)$ & $ (1, 1)$ & $\dots$ & 
 $(1,1) $ & $1$ & $ 0 $
\\
$d $ & $ =$ & $( 0, 2)$ &$(0, 2)$ & $\dots$ & 
 $(0,2)$ & $(0,2)$ & $ (0,2)$ & $\dots$ & 
 $(0,2) $ & $0 $ &$ 1 $
\end{tabular}
  \end{equation*}    

The element $c$ witnesses that 
$(a,d) \in  \beta \circ \gamma $. Moreover,
$ a \mathrel { \alpha  } d$, hence 
$(a,d) \in  \alpha (\beta \circ \gamma ) $.
On the other hand, 
the only other element 
$ \alpha \beta $-connected  to $a$ is
\begin{equation*}
\begin{tabular}{lccccccccccc}
$f_1 $ & $ =$ & $( 1, 0)$ &$(2,0)$ & $\dots$ & 
 $(2,0)$ & $( 2, 0)$ & $ (2, 0)$ & $\dots$ & 
 $(2,0) $ & $2 $ &$ 1 $
\end{tabular}
  \end{equation*}    
due to the rule (b3) in the formation of $P$.
Due to the definition of $\gamma$,
the only other element 
$  \gamma  $-connected  to $f_1$ 
and with last component $1$ is
\begin{equation*}
\begin{tabular}{lccccccccccc}
$f_2 $ & $ =$ & $( 0, 0)$ &$(2,0)$ & $\dots$ & 
 $(2,0)$ & $( 2, 0)$ & $ (2, 0)$ & $\dots$ & 
 $(2,0) $ & $2 $ &$ 1 $
\end{tabular}
  \end{equation*}

Continuing this way, the only possibility
to go from $a$ to  $d$ 
 through an $ \alpha \beta $-or-$ \alpha \gamma $-chain 
is to consider \emph{all} the elements

\begin{equation*}
\begin{tabular}{rccccccccccc}
$a $ & $ =$ & $( 2, 0)$ &$(2,0)$ & $\dots$ & 
 $(2,0)$ & $( 2, 0)$ & $ (2, 0)$ & $\dots$ & 
 $(2,0) $ & $2 $ &$ 1 $
\\
$f_1 $ & $ =$ & $( 1, 0)$ &$(2,0)$ & $\dots$ & 
 $(2,0)$ & $( 2, 0)$ & $ (2, 0)$ & $\dots$ & 
 $(2,0) $ & $2 $ &$ 1 $
\\
$f_2 $ & $ =$ & $( 0, 0)$ &$(2,0)$ & $\dots$ & 
 $(2,0)$ & $( 2, 0)$ & $ (2, 0)$ & $\dots$ & 
 $(2,0) $ & $2 $ &$ 1 $
\\
$f_3 $ & $ =$ & $( 0, 0)$ &$(1,0)$ & $\dots$ & 
 $(2,0)$ & $( 2, 0)$ & $ (2, 0)$ & $\dots$ & 
 $(2,0) $ & $2 $ &$ 1 $
\\
 & &  & & & 
  & \dots &  &  & 
  &  &
\\
$f_{m-4} $ & $ =$ & $(0,0)$ &$(0,0)$ & $\dots$ & 
 $(0,0)$ & $( 0, 0)$ & $ (0, 0)$ & $\dots$ & 
 $(1,0) $ & $ 2 $ &$ 1 $
\\
$f_{m-3}$ & $ =$ & $(0,0)$ &$(0,0)$ & $\dots$ & 
 $(0,0)$ & $( 0, 0)$ & $ (0, 0)$ & $\dots$ & 
 $(0,0) $ & $ 2 $ &$ 1 $
\\
$f_{m-2}$ & $ =$ & $(0,0)$ &$(0,0)$ & $\dots$ & 
 $(0,0)$ & $( 0, 0)$ & $ (0, 0)$ & $\dots$ & 
 $(0,0) $ & $ 1 $ &$ 1 $
\\
$f_{m-1}$ & $ =$ & $(0,0)$ &$(0,0)$ & $\dots$ & 
 $(0,0)$ & $( 0, 0)$ & $ (0, 0)$ & $\dots$ & 
 $(0,0) $ & $ 0 $ &$ 1 $
\\
$f_m $ & $ =$ & $(0,0)$ &$(0,0)$ & $\dots$ & 
 $(0,0)$ & $( 0, 0)$ & $ (0, 0)$ & $\dots$ & 
 $(0,1) $ & $ 0 $ &$ 1 $
\\
 & &  & & & 
  & \dots &  &  & 
  &  &
\\
$f_{2m-7} $ & $ =$ & $( 0, 0)$ &$(0, 1)$ & $\dots$ & 
 $(0,2)$ & $(0,2)$ & $ (0,2)$ & $\dots$ & 
 $(0,2) $ & $0 $ &$ 1 $
\\
$f_{2m-6} $ & $ =$ & $( 0, 0)$ &$(0, 2)$ & $\dots$ & 
 $(0,2)$ & $(0,2)$ & $ (0,2)$ & $\dots$ & 
 $(0,2) $ & $0 $ &$ 1 $
\\
$f_{2m-5}$ & $ =$ & $( 0, 1)$ &$(0, 2)$ & $\dots$ & 
 $(0,2)$ & $(0,2)$ & $ (0,2)$ & $\dots$ & 
 $(0,2) $ & $0 $ &$ 1 $
\\
$d $ & $ =$ & $( 0, 2)$ &$(0, 2)$ & $\dots$ & 
 $(0,2)$ & $(0,2)$ & $ (0,2)$ & $\dots$ & 
 $(0,2) $ & $0 $ &$ 1 $
\end{tabular}
  \end{equation*}
in the case $m$ odd, while in the case $m$ even 
the penultimate column should be deleted and
the 
above ``middle'' block is replaced by
\begin{equation*}
\begin{tabular}{rcccccccccc}
 & &  & & & 
  & \dots &  &  & 
  &  
\\
$f_{m-3} $ & $ =$ & $(0,0)$ &$(0,0)$ & $\dots$ & 
 $(0,0)$ & $( 0, 0)$ & $ (0, 0)$ & $\dots$ & 
 $(1,0) $  &$ 1 $
\\
$f_{m-2}$ & $ =$ & $(0,0)$ &$(0,0)$ & $\dots$ & 
 $(0,0)$ & $( 0, 0)$ & $ (0, 0)$ & $\dots$ & 
 $(0,0) $  &$ 1 $
\\
$f_{m-1} $ & $ =$ & $(0,0)$ &$(0,0)$ & $\dots$ & 
 $(0,0)$ & $( 0, 0)$ & $ (0, 0)$ & $\dots$ & 
 $(0,1) $  &$ 1 $
\\
 & &  & & & 
  & \dots &  &  & 
  &  
\end{tabular}
  \end{equation*}
   
Since we need  to consider all the above elements,
we get
 $(a,d)  \notin \alpha \beta 
\circ \alpha \gamma \circ {\stackrel{2m-5}{\dots}}  $, hence
 $\alpha (\beta \circ \gamma )  \not \subseteq  \alpha \beta 
\circ \alpha \gamma \circ {\stackrel{2m-5}{\dots}}   $, that is, 
$\mathbf B$ does not belong to a $2m{-}5$-distributive variety.   
Recall that
$ \alpha \beta \circ \alpha  \gamma \circ {\stackrel{k}{\dots}}  $
denotes the relation $ \alpha \beta \circ \alpha \gamma
 \circ \alpha \beta \circ \dots$
with $k-1$ occurrences of  $ \circ $
and that $R^k$ is $R \circ R \circ {\stackrel{k}{\dots}}  $  

As implicit in the proof of \ref{mmm}, we see that $a$
is $ \alpha \gamma $-connected only to itself, hence
we also get   
$\alpha (\beta \circ \gamma )  \not \subseteq  \alpha \gamma 
\circ \alpha \beta  \circ {\stackrel{2m-4}{\dots}}    $, that is, 
$\mathbf B$ does not belong to a $2m{-}4$-alvin variety.

In the above arguments we have considered identities involving 
$\alpha \gamma \circ \alpha \beta  \circ \dots$ on the right
only for simplicity. 
While $a$ is $ \beta $-connected to further elements of 
$B$, since we can consider elements with $0$ as the last coordinate,
on the other hand, $f_1$ is the only other element such that 
$ (a, f_1) \in  \alpha ( \gamma  \circ \beta  ) $.
Continuing the same way,
the only elements $h$ such that  
$ (f_1, h) \in  \alpha ( \gamma  \circ \beta  ) $
are $a$,  $f_2$ and  $f_3$.
Of course, it is no use to turn back to $a$, and the ``fastest way to $d$''
uses $f_3$.
Going on, we see that     
 $(a,d)  \notin (\alpha (  \gamma \circ \beta )) ^{m-2} $,
 hence
 $\alpha (\beta \circ \gamma )  \not \subseteq 
 (\alpha ( \gamma \circ \beta )) ^{m-2}$,
that is, clause (4) in Theorem \ref{mmm}.

Dealing with larger   even $q$ presents no significant difference,
while if $q$ is odd it is enough to modify 
the definitions displayed in \eqref{gath}:
all the pairs in the definitions of $\beta$ and $\gamma$  
should be, respectively,
 $(\beta^*_{q+1} \times \beta  ^*_{q+1} )$    and
 $( \gamma  ^*_{q+1} \times \gamma   ^*_{q+1} )$.
Then it is more convenient to deal with
identity \eqref{ciodd} in Lemma \ref{cilem}.  
For example, the first elements in the shortest chain from
$a$ to  $d$ in the case $q=3$ are
\begin{equation*}
\begin{tabular}{rccccccccccc}
$a $ & $ =$ & $( 3, 0)$ &$(3,0)$ & $\dots$ & 
 $(3,0)$ & $( 3, 0)$ & $ (3, 0)$ & $\dots$ & 
 $(3,0) $ & $3 $ &$ 1 $
\\
$f_1 $ & $ =$ & $( 2, 0)$ &$(3,0)$ & $\dots$ & 
 $(3,0)$ & $( 3, 0)$ & $ (3, 0)$ & $\dots$ & 
 $(3,0) $ & $3 $ &$ 1 $
\\
$f_2 $ & $ =$ & $( 1, 0)$ &$(3,0)$ & $\dots$ & 
 $(3,0)$ & $( 3, 0)$ & $ (3, 0)$ & $\dots$ & 
 $(3,0) $ & $3 $ &$ 1 $
\\
$f_3 $ & $ =$ & $( 0, 0)$ &$(2,0)$ & $\dots$ & 
 $(3,0)$ & $( 3, 0)$ & $ (3, 0)$ & $\dots$ & 
 $(3,0) $ & $3 $ &$ 1 $
\\
 & &  & & & 
  & \dots &  &  & 
  & 
\end{tabular}
  \end{equation*}

Notice that, as in the case $q$ even, 
$a \mathrel { \beta  } f_1 $,
but there is no other element $\alpha \gamma $-connected
to $a$. On the other hand, in the case $q$ odd
we are sometimes able to move two components at a time, 
as is the case for $f_2$ and  $f_3$ above.   
   
In conclusion, the above arguments show that the following proposition 
holds. Recall that the definitions of $\mathbf P = \mathbf P (m,q)$ and 
 $\mathbf B =  \mathbf B (m,q)$ 
depend on $m$ and $q$, though sometimes we have not  explicitly
indicated the dependence in the above arguments.

\begin{proposition} \labbel{b}
For every $m \geq 3$ and $q \geq 2$, the algebra 
$\mathbf B(m,q)$, as constructed above, has an
$m$-ary near-unanimity term. Identity  
  \eqref{cidd} fails in $\mathbf B(m,q)$. 
 \end{proposition}

\section{Further remarks} \labbel{fur} 

It is well-known that, for every $m \geq 4$,
there is a variety with an $m$-ary near-unanimity term   
and without an $m{-}1$-ary near-unanimity term.
See, e.~g.,  \cite[Lemma 3.4]{S} or Example \ref{ex} here.
The variety $\mathcal {N}_m$ introduced in Definition \ref{vars}
  furnishes  another counterexample, as
we shall show in the next corollary.
In addition, the counterexamples presented here
  have a 
symmetric $m$-ary near-unanimity term and
are locally finite.
Notice that the variety denoted by $\mathcal {N}_m$
in \cite{S} is distinct from the variety denoted by   $\mathcal {N}_m$
here. Also, the indices are shifted by $1$ in most definitions, with respect to 
\cite{S}. 

For every $n \geq 2$ there are known examples
of $n$-distributive not $n{-}1$-distributive varieties,
e.~g. \cite{FV,day} and further references there.
For $n$ even, the variety
$\mathcal {N}_m$ provides another example (with
$ m=\frac{n+4}{2}$).
Corresponding examples appear in  \cite{day}
regarding $n$-modularity; again,
$\mathcal {N}_m$ provides further counterexamples.

\begin{corollary} \labbel{n-1}
If $m \geq 4$, then $\mathcal {N}_m$ has 
a symmetric $m$-ary near-unanimity term   but
no $m{-}1$-ary near-unanimity term (symmetric or not).

The variety  $\mathcal {N}_m$
is $2m{-}4$-distributive but not
$2m{-}5$-distributive,
$2m{-}3$-modular but not
$2m{-}4$-modular.
 \end{corollary}

 \begin{proof}
We have proved in 
Theorem \ref{mmm}
that $\mathcal {N}_m$ has 
an $m$-ary symmetric near-unanimity term. 
The arguments in Example \ref{ex}
show that the variety generated by 
the algebra $\mathbf N ^{2,m} $ from Definition \ref{vars}
has not  an $m{-}1$-ary near-unanimity term.
Since $\mathbf N ^{2,m} $ is one among the generators
of $\mathcal {N}_m$, then $\mathcal {N}_m$ 
has not  an $m{-}1$-ary near-unanimity term.

The second statement is immediate from
Theorems  \ref{mnun}
and \ref{mmm}(1)(3).

From Theorems  \ref{mnun}
and \ref{mmm}(1)(3) we can 
also obtain
another proof 
 that $\mathcal {N}_m$ has not
an $m{-}1$-ary near-unanimity term.
If, by contradiction,  such a term exists, then  $\mathcal {N}_m$
is $2m{-}6$-distributive.   
However, $\mathcal {N}_m$
is not even $2m{-}5$-distributive, absurd.
A similar proof is obtained by dealing with
the modularity levels.
Remark \ref{mult}  could  be used to get still another
proof. 
 \end{proof} 

\begin{remark} \labbel{n'}
Recall the definition of 
$\mathbf L ^{{\rm nu}, j,m} $ from Definition \ref{vars}.
 
For any given $m \geq 3$,
let $\ell = \frac{m+1}{2} $
if $m$ is odd, 
  $\ell = \frac{m}{2} $
if $m$ is even
and 
 let $\mathcal {N}'_m$ be the variety 
generated by all the  lattice
reducts 
$\mathbf L ^{{\rm nu}, j,m} $,
with $\mathbf L$ an arbitrary lattice, 
$2 \leq j \leq \ell $. 

Theorem \ref{mmm}
and Corollary \ref{n-1}
hold for $\mathcal {N}'_m$, as well,
with the only exception that $\mathcal {N}'_m$
is not locally finite.
 \end{remark}

In another classical paper 
Mitschke \cite{mim} 
showed that the variety $\mathcal I$ of implication algebras
is congruence $3$-distributive, $3$-permutable, not  
$2$-distributive and not permutable.
Then in \cite{Mi} she proved that for no $m$  
the variety $\mathcal I$ has an $m$-ary  near-unanimity term.
Another proof can be found in \cite[Remarks 2.2(a)(b)]{B}.

In \cite[Section 5]{ia} we expanded   $\mathcal I$
by adding a $4$-ary near-unanimity term 
in such a way that the distributive and permutable levels
remain unchanged.
 Combining the arguments from \cite{mim,Mi,ia} and 
from Example \ref{ex}, we  get a 
variety sharing the same levels of $\mathcal I$, with an
$m+1$-ary  near-unanimity term
  but  without an $m$-ary  near-unanimity term.
See the next proposition.

This shows that $3$-permutability 
has no effect on the integers $m$
for which an $m$-ary  near-unanimity can exist.  
Notice that, on the other hand, 
it is immediate from the characterizations
in \eqref{cacu} that a congruence distributive
(modular) $n$-permutable variety is $n$-distributive
($n$-modular).   
Notice also that a congruence permutable variety
with a near-unanimity term is congruence distributive,
hence, by permutability, $2$-distributive,
and this means the existence of a majority term, namely, a $3$-ary
near-unanimity term. Thus, in contrast with
$3$-permutability, the stronger notion of 
permutability does trivialize
 the sets of integers $m$
for which an $m$-ary  near-unanimity  exists.

The operations of a Boolean algebra  
shall be  denoted by $+$, $\cdot$ and $'$.
The variety $\mathcal I$ of (dual) \emph{implication algebras}
is the variety generated by term-reducts of Boolean algebras
in which $i(x,y) = xy'$ is the only basic operation.
Let $f(x,y,z) $ be the Boolean term $  x(y'+z)$.
The variety 
$\mathcal I^-$
is the variety generated by reducts of Boolean algebras 
having  $f$  as the only basic operation.
If $m \geq 3$, we let $\mathcal I_m$, resp. $\mathcal I^-_m$,
be the varieties generated by reducts of Boolean algebra 
with two basic operations  corresponding to $i$ and  
$u _{2,m} $, resp.,  $f$ and  
$u _{2,m} $. Recall the definition of   
$u _{2,m} $  from Definition \ref{latterms}.

\begin{proposition} \labbel{MM}
If $m \geq 4$, then  both $\mathcal I_{m}$ and 
 $\mathcal I^{-}_m$
are  $3$-distributive, congruence $3$-permutable,
not congruence permutable, not
$2$-distributive, have an $m$-ary near-unanimity term  
but not an $m{-}1$-ary near-unanimity term. 
 \end{proposition}

  \begin{proof} 
Since 
$\mathcal I$ and 
 $\mathcal I^{-}$ are
$3$-distributive and congruence $3$-permutable
 \cite{mim}, then so are their expansions $\mathcal I_{m}$ and 
 $\mathcal I^{-}_m$.
To show that $2$-distributivity and permutability fail, consider
the reduct  $\mathbf A$ of $\mathbf 2^3$, where $\mathbf 2$
is the two-elements Boolean algebra.
Since $m \geq 4$, then  $2^3 \setminus (1,1,1)$
is closed   under $u _{2,m} $, compare an argument in Example \ref{ex}. 
Moreover, 
$2^3 \setminus (1,1,1)$ is closed also under 
$i$ (or $f$), thus it is the universe for a subalgebra of 
 $\mathbf A$. The original argument in \cite{mim}
(credited in that form to the referee)
now shows that $2$-distributivity fails,
hence also permutability fails.
The argument is recalled also in the proof of  \cite[Proposition 5.1]{ia}. 

To show that neither $\mathcal I_{m}$ nor 
 $\mathcal I^{-}_m$
have an $m-1$-ary near-unanimity term
argue as in Example \ref{ex}, considering 
the reduct of  
$\mathbf 2^{m-1} \setminus (1,1, \dots , 1) $.
\end{proof}

\begin{remark} \labbel{BaKo}   
For $m \geq 3$, 
an  $m$-ary near-unanimity term 
implies the existence of a sequence
$t_1, \dots, t_{m-2}$ of \emph{directed J{\'o}nsson terms},
i. e., terms satisfying
\begin{align*}  
x&=t_1(x,x,z), 
&& \quad  t_{m-2}(x,z,z)=z, 
\\
t_{i}(x,z,z)&=  t_{i+1}(x,x,z), &&\text{for $1 \leq i <m-2$, and  }
\\
 x&= t_{i}(x,y,x),  && \text{for $1 \leq i \leq m-2$.}
   \end{align*}    
See  Barto and Kozik \cite[Section 5.3.1]{BK}.
Directed J{\'o}nsson terms provide another characterization
of congruence distributivity \cite{KKMM}.

The mentioned observation from 
\cite[Section 5.3.1]{BK} is optimal:
the variety $\mathcal {N}_m$
fails to have
$t_1, \dots, t_{m-3}$ directed J{\'o}nsson terms, since
otherwise $\mathcal {N}_m$ would be 
$2m{-}6$-distributive, by \cite[Observation 1.2]{KKMM}, 
thus contradicting   
Theorem \ref{mmm}(1). 
Notice that the counting conventions
in \cite{BK,KKMM} are sometimes different from 
the conventions adopted in the present note. 
 \end{remark}

\begin{remark} \labbel{nat}
Theorem \ref{mnun}
suggests that, for every 
  $m \geq 3$, there should be  some 
condition $C_m$ strictly between the strength of 
an   $m$-ary and of an $m{+}1$-ary
near-unanimity term and such that 
$C_m$ implies $2m{-}3$-distributivity but does not imply
$2m{-}4$-\brfrt distributivity.
Of course,

($\diamondsuit_m$) there is an $m{+}1$-ary
near-unanimity term + $2m{-}3$-distributivity

\noindent
is possibly such a condition, but it looks quite artificial.
One should check that 
($\diamondsuit_m$) 
does not imply  
$2m{-}4$-distributivity,
a fact which can be probably obtained 
by combining the present methods with \cite{day}.
For $m=3$, a $3$-distributive not $2$-distributive
variety with a $4$-ary near-unanimity term does indeed 
exist, see \cite[Proposition 5.1]{ia}
or Proposition \ref{MM} above in the case $m=4$.        

A possibly more natural condition is presented in the next definition.
 \end{remark}  

\begin{definition} \labbel{mez}
If $m \geq 3$, an \emph{$m \frac{1}{2} $-near-unanimity term}
is an $m{+}2$-ary term $u $ such that the following equations hold.
\begin{align} 
\labbel{a1}   
&u(z,z, x, x, \dots, x) = x,
\\
\labbel{a3}
&u(x, \dots, x,\underset{i}{z},x, \dots, x) = x, && \text{for $2 \leq i \leq m+2$,} 
\\ \labbel{a2}
&u(x, x, x, z, z, \dots, z) = u(x, z, z, z, z, \dots, z).
 \end{align}   

The terminology comes from the fact that 
if $u$ is an $m \frac{1}{2} $-near-unanimity term,
then the $m{+}1$-ary term $v$ defined by 
\begin{equation*}\labbel{v}
v(x_1, x_2, x_3, \dots , x_{m+1})=
 u(x_1, x_1, x_2, x_3, \dots , x_{m+1})
 \end{equation*}    
is a near-unanimity term, by \eqref{a1} and \eqref{a3}.
On the other hand,   
if $w$ is an $m$-ary near-unanimity term,
then by adding two initial dummy varables,
$w$ becomes    an
$m \frac{1}{2} $-near-unanimity term.
 \end{definition}

\begin{proposition} \labbel{propmez}
Let $m \geq 3$.
  \begin{enumerate}   
 \item 
If some variety $\mathcal V$ has an
$m \frac{1}{2} $-near-unanimity term, then 
$\mathcal V$ is $2m{-}3$-distributive.
\item 
There is a variety  $\mathcal V$ 
with an $m{+}1$-ary near-unanimity term 
but without an
$m \frac{1}{2} $-near-unanimity term.
\item
There is a variety  $\mathcal V$ 
with an 
$m \frac{1}{2} $-near-unanimity term
but without an $m$-ary near-unanimity term. 
  \end{enumerate}  
 \end{proposition} 

Details for the proof of Proposition \ref{propmez} 
 shall be presented  elsewhere.
We just notice that, granted clause (1)  in Proposition \ref{propmez},
then, by
 Theorem \ref{mmm}(1), the variety $\mathcal N_{m+1}$
furnishes an example for \ref{propmez}(2).

\begin{problem} \labbel{diss}
Study ``dissent'' terms in the following sense.
A term $u$ of arity $\geq 3$ is a \emph{lone-dissent term}  
if  all the equations of the form
\begin{equation*} 
u(x,x, \dots, x,y,x, \dots, x,  x) =y
\end{equation*}    
are satisfied, with just one occurrence of $y$
in any possible position. 

A ternary lone-dissent term is a \emph{minority term},
see \cite{KOVZ} for a detailed study of minority terms. 
 It is easy to see that a variety of abelian groups
has an $m{+}1$-ary lone-dissent term if and only if 
its exponent divides $m$.
In particular, contrary to the case of near-unanimity terms,
the existence of an $m$-ary lone-dissent term
does not imply an $m{+1}$-ary lone-dissent term.
However, there are some  positive results.
Some simple facts are stated in the next proposition.

It is probably also interesting to study 
dissent-unanimity terms in the following sense.
If $m\geq 3$, a $2m$-ary term $u$ is a    
\emph{dissent-unanimity term}  
if  all the equations of the form
\begin{equation*} 
u(x,x, \dots, x,\underset{i}{y},x, \dots, x,  x;
 y,y, \dots, y,\underset{m+i}{z},y, \dots , y, y) =y
\end{equation*}    
are satisfied for all $i$, $1 \leq i \leq m$,
where the semicolon separates the first 
$m$ arguments of $u$ with the last $m$ arguments.  
The case $m=3$ has been dealt with in \cite{KOVZ}.   
\end{problem}   
 
\begin{proposition} \labbel{propdiss}
Let $m, n \geq 2$.   
\begin{enumerate} 
   \item 
If some variety $\mathcal V$ has an
$m{+}1$-ary lone-dissent term, then, for all $k \geq 1$,
$\mathcal V$ has a $km{+}1$-ary lone-dissent term.
\item
More generally,
if some variety $\mathcal V$ has both an
$m{+}1$-ary and an $n{+}1$-ary lone-dissent term,
then $\mathcal V$ 
has 
an $m{+}n{+}1$-ary lone-dissent term.
\item 
If some variety $\mathcal V$ has an
$m{+}1$-ary lone-dissent term, then $\mathcal V$ has a Maltsev
term, hence $\mathcal V$ is congruence permutable.
\item
If some variety $\mathcal V$ has both an
$m{+}1$-ary and an $m{+}2$-ary lone-dissent term,
then $\mathcal V$ 
is an arithmetical variety. 
  \end{enumerate}
 \end{proposition}

  \begin{proof}  
 (1) If $d$ is an $m{+}1$-ary lone-dissent term,
then $d(d(\hyphe,\hyphe, \dots ), \hyphe, \dots)$
is a $2m{+}1$-ary lone-dissent term,   
$d(d(d(\hyphe,\hyphe, \dots ),\hyphe, \dots ), \hyphe, \dots)$
is a $3m{+}1$-ary lone-dissent term and so on.

(2) If $d$, $e$  are, respectively, 
 an  $m{+}1$-ary and an $n{+}1$-ary lone-dissent term,
then $d(e(\hyphe,\hyphe, \dots ), \hyphe, \dots)$
is an $m{+}n{+}1$-ary lone-dissent term.

(3)  If $d$ is an $m{+}1$-ary lone-dissent term,
then $t(x,y,z)=d(x,y,y, \dots, y, z)$ is a Maltsev term.

(4) If $d$, $e$  are an  $m{+}1$-ary and an $m{+}2$-ary lone-dissent term,
then 
\begin{align*}
 t(&x_1, \dots , x_{m+2})=
\\
e(&d(x_1, x_2, \dots, x_{m}, x_{m+1}), 
d(x_1, x_2, \dots, x_{m}, x_{m+2}),
\\
&d(x_1, x_2, \dots, x_{m-1},x_{m+1}, x_{m+2}), \dots,
d(x_1, x_2, x_4, \dots, x_{m+1}, x_{m+2}),
\\ 
&d(x_1, x_3,  \dots, x_{m+1}, x_{m+2}),
d(x_2, x_3, \dots, x_{m+1}, x_{m+2}))
\end{align*}
is an $m{+}1$-ary near-unanimity term, hence
$\mathcal V$ is congruence distributive. By (3)
$\mathcal V$ is congruence permutable, hence $\mathcal V$ 
is arithmetical.  
\end{proof}

\begin{corollary} \labbel{cordiss}
If $m,n \geq 2$, $m$ and $n$ are coprime and 
some variety  $\mathcal V$ has both an
$m{+}1$-ary and an $n{+}1$-ary lone-dissent term,
then $\mathcal V$ 
is arithmetical. 
 \end{corollary} 

\begin{proof} 
 Since $m$ and $n$ are coprime,
then the Diophantine equation $km+hn=1$
has a solution with $k,h \in \mathbb Z$. 
Changing a sign, we have either 
$km=hn+1$
or 
$hn=km+1$ with $k,h \in \mathbb N \setminus \{ 0 \} $.
Applying \ref{propdiss}(1) twice, $\mathcal V$ 
has both a $km{+}1$- and 
a $hn{+}1$-lone-dissent term.
Since $hn$ and $km$ differ by $1$,
then $\mathcal V$ is arithmetical by \ref{propdiss}(4).    
\end{proof}  

As a final remark, we notice that, for $m \geq 3$,
the existence of an $m$-ary lone-dissent term
implies the existence of a sequence   
$t_1, \dots, t_{m-2}$ of \emph{directed minority terms},
i. e., terms satisfying
\begin{align*} \labbel{dirdis} 
y&=t_1(x,x,y), 
&& \quad  t_{m-2}(x,y,y)=x, 
\\
t_{i}(x,y,y)&=  t_{i+1}(x,x,y), &&\text{for $1 \leq i <m-2$, and  }
\\
 y&= t_{i}(x,y,x),  && \text{for $1 \leq i \leq m-2$.}
   \end{align*}    
The proof presents no variation with respect to  \cite[Section 5.3.1]{BK}.
It is probably interesting to study this and similar conditions.
See also \cite[Remark 8.19]{day}
for further comments.

\smallskip 

{\scriptsize
The author 
 considers  highly  inappropriate 
and strongly discourages the use 
of indicators extracted from the following list 
(even in aggregate forms in combination with similar lists)
 in decisions about individuals,
 attributions of funds, selections or evaluations of research projects.
\par
}


\begin{thebibliography}{10}

\bibitem{Bk} Baker, K.A.:
Congruence-distributive polynomial reducts of lattices.
Algebra Universalis
\textbf{9}, 142--145
(1979)

\bibitem{BP}
Baker, K.A., Pixley, A.F.: Polynomial interpolation and the Chinese remainder
theorem for algebraic systems. Math. Z. \textbf{143}, 165--174 (1975) 

\bibitem{Ba} Barto, L.: Finitely related algebras in congruence distributive varieties have near unanimity terms. Canad. J. Math. \textbf{65}, 3--21 (2013) 

\bibitem{BK}
Barto, L., Kozik, M.: `Absorption in universal algebra and CSPV', in: The Constraint Satisfaction
Problem: Complexity and Approximability, Dagstuhl Follow-Ups, 7 (Schloss 
Dagstuhl--Leibniz
Zentrum f\"ur Informatik, Wadern, 2017), 45--77

\bibitem{BIM} Berman, J., Idziak, P., Markovi\'c, P., McKenzie, R., Valeriote, M., Willard, R.:
Varieties with few subalgebras of powers. Trans. Amer. Math. Soc. \textbf{362},
1445--1473 (2010) 


 \bibitem{CCV} Campanella, M., Conley, S., Valeriote, M.:
{Preserving near unanimity terms under products.}
Algebra Universalis \textbf{76}, 293--300  (2016) 
 

\bibitem{D} 
 Day, A.: {A characterization of modularity for congruence lattices of algebras.}
 Canad. Math.
Bull. \textbf{12}, 167--173 (1969)


\bibitem{FV} Freese, R., Valeriote, M.A.:
 On the complexity of some Maltsev conditions. Int.
J. Algebra Comput. \textbf{19}, 41--77 (2009)

\bibitem
{JD} J{\'o}nsson, B.: {Algebras whose congruence lattices
   are distributive.} Math. Scand. \textbf{21}, 110--121  (1967)

\bibitem{KP} Kaarli, K., Pixley,  A.F.:  Polynomial Completeness 
in Algebraic Systems (Chapman \& Hall/CRC,
Boca Raton, FL, 2001)

\bibitem{KKMM}
Kazda, A., Kozik, M., McKenzie, R., Moore, M.: Absorption and directed
J\'onsson terms. In: Czelakowski, J. (ed.) Don Pigozzi on Abstract Algebraic
Logic, Universal Algebra, and Computer Science, Outstanding Contributions to
Logic \textbf{16}, pp. 203--220. Springer, Cham (2018) 

\bibitem{KOVZ} 
  Kazda, A., Opr\v sal, J.,  Valeriote,  M., Zhuk, D.:
Deciding the existence of minority terms.
Canad. Math.
Bull. \textbf{63}, 577--591 (2020). 


\bibitem{ia} Lipparini, P.:
 Relation identities in 3-distributive varieties,
Algebra Universalis
\textbf{80},
Paper No. 55, 20
(2019)


\bibitem
{B} Lipparini, P.: 
The distributivity spectrum of Baker's variety, 
J. Aust. Math. Soc. \textbf{110},
119--144 (2021)

\bibitem{day}
Lipparini, P.: {Day's Theorem is sharp for $n$ even.} arXiv:1902.05995,
1--60 (2019/2021)

\bibitem
{MMT}  McKenzie, R.N.,  McNulty, G.F., 
Taylor, W. F.: {Algebras, Lattices, Varieties. Vol. I}
 (Wadsworth \& Brooks/Cole
   Advanced Books \& Software, 1987),
corrected reprint with additional bibliography (AMS Chelsea Publishing/American
Mathematical Society, 2018)


\bibitem{mim} Mitschke, A.: Implication Algebras are
 $3$-Permutable and $3$-Distributive. Algebra
Universalis \textbf{1}, 182--186 (1971)


\bibitem{Mi} 
Mitschke, A.:
{Near unanimity identities and congruence distributivity in equational classes.}
 Algebra Universalis \textbf{8},  29--32  (1978)


\bibitem{S}
Sequeira, L.:
{Near-unanimity is decomposable.}
Algebra Universalis  \textbf{50},  157--164 (2003) 

\bibitem{Ts} Tschantz, S.T.: More conditions equivalent to congruence modularity.
 In: Universal Algebra and Lattice Theory, 270--282, Lecture Notes in Math. 
\textbf{1149} (1985)


\end{thebibliography}
\end{document}